\RequirePackage{fix-cm}
\documentclass[twocolumn]{svjour3}          
\usepackage{algorithm}
\usepackage{algorithmic}
\usepackage{amssymb}
\usepackage{graphicx}
\usepackage{diagbox}
\usepackage{amsmath}
\usepackage[english]{babel}
\usepackage[T1]{fontenc}
\usepackage[utf8]{inputenc}
\usepackage{placeins}
\usepackage{url}
\usepackage{tikzpagenodes}
\usepackage{tikz}
\usepackage{lineno}
\usepackage{commath}
\usetikzlibrary{arrows.meta}

\usetikzlibrary{arrows,matrix,positioning}

\DeclareMathOperator*{\argmin}{arg\,min}

\smartqed  

\begin{document}

\title{Model Order Reduction for\newline Efficient Descriptor-Based Shape Analysis
}

\titlerunning{Model Order Reduction for Efficient Descriptor-Based Shape Analysis}     

\author{Martin Bähr \and Michael Breuß \and Robert Dachsel }

\authorrunning{M. Bähr et al.} 

\institute{Brandenburg Technical University, Institute for Mathematics,\\
Platz der Deutschen Einheit 1, 03046 Cottbus, Germany, \\
\email{$\lbrace$bähr,breuss,dachsel$\rbrace$@b-tu.de}}


\date{Received: date / Accepted: date}

\maketitle
\begin{abstract}
In order to investigate correspondences between 3D shapes,
many methods rely on a feature descriptor which is invariant under
almost isometric transformations. An interesting class of models for
such descriptors relies on partial differential equations (PDEs) 
based on the Laplace-Beltrami operator for constructing intrinsic shape signatures. 
In order to conduct the construction, not only a variety of PDEs but also 
several ways to solve them have been considered in previous works. In particular, 
spectral methods have been used derived from the series expansion of analytic solutions
of the PDEs, and alternatively numerical integration schemes have been proposed.

In this paper we show how to define a computational framework by model order reduction
(MOR) that yields efficient PDE integration
and much more accurate shape signatures as in previous works.
Within the construction of our framework we introduce some technical
novelties that contribute to these advances, and in doing this we present some improvements for
virtually all considered methods. As part of the main contributions, we show
for the first time an extensive and detailed comparison between the spectral and integration techniques,
which is possible by the advances documented in this paper. We also propose here to employ soft correspondences
in the context of the MOR methods which turns out to be highly beneficial with this approach.
\keywords{Shape Matching, Shape Analysis, Partial Differential Equations, Heat Equation, Wave Equation,
Model Order Reduction, Soft Correspondence}
\end{abstract}

\section{Introduction}

The investigation of correspondences between 3D shapes is a fundamental problem in computer
vision, graphics and pattern recognition. It has a wide variety of potential applications, 
including e.g.\ shape comparison or texture transfer, see e.g.\ \cite{vanKaick-etal-2011} for
some discussion. The basic task of finding shape correspondences
is to identify a relation between elements of two or more shapes, and a challenging setting for
this is concerned with non-rigid shapes that are assumed to be just almost isometric,
compare for instance \cite{TOSCA}.

One of the possible strategies to find pointwise correspondences is to construct a feature descriptor, 
or shape signature, which characterises geometry around the points that define the surface 
of a given shape. An interesting type of such approaches relies on the Laplace-Beltrami operator 
which enables to describe intrinsic geometric properties of a shapes' surface \cite{Levy-2006,Reuter-etal-2009,R}.
To this end, several partial differential equations (PDEs) have been proposed for the construction of shape signatures
that rely on the Laplace-Beltrami operator as a crucial component. The {\em heat kernel signature} (HKS)
introduced in \cite{SOG} is based on the heat equation, the {\em wave kernel signature} (WKS)
\cite{ASC} is constructed using the Schr\"odinger equation, and recently the classic
wave equation \cite{DBH-2017} has been applied. 

The mentioned kernel-based signatures (HKS and WKS) are computed using the kernels of analytic solutions of the 
underlying PDEs, resulting in infinite series expansions that define the complete solutions in time at 
a given surface point.
Let us note that it is not immediately evident when to truncate the series for practical purposes, yet some 
strategies have been given in the literature \cite{TOSCA}; see also \cite{DBH-2018} for a recent,
related investigation. As an alternative approach to the kernel-based methods, 
the full numerical integration of the PDEs has been considered \cite{DBH-2016}, and the
abovementioned shape signature defined via the classic wave equation \cite{DBH-2017} has been realised in
this setting. Comparing these two approaches,
the numerical integration as reported in \cite{DBH-2017,DBH-2016} may yield a more accurate correspondence, 
while the computation of the numerical signatures has been much more time-consuming than the kernel-based signatures. 
In order to improve the performance of the integration approach, a {\em model order reduction} (MOR) technique
has been proposed in \cite{BDB-2018}, showing some promising first results.

\paragraph{Our Contribution}

The goal of this paper is to document several substantial advances of the integration based approach. Based on these
improvements, we also give here a much more detailed comparison to kernel-based methods than it was possible before.
More concretely, let us mention the following main points, clarifying at the same time relation to previous literature:
\par
    {\em (i)} Employing the MOR technique introduced in \cite{BDB-2018} for the shape matching purpose,
    we show that it pays off to extend the considered time scale in the numerical MOR signatures.
    The {\em large time scales} we employ here for the first time with the MOR signature yields a notable
    improvement in quality of results.
\par
    {\em (ii)} We introduce a {\em numerical similarity transform} within both the MOR and kernel-based signatures
    which yields a much more stable computation as well as in most cases results of much higher accuracy
    for {\em all} considered techniques.
\par
    {\em (iii)} Since the MOR method is highly efficient, it is now possible to evaluate much larger datasets
    conveniently. We present here for the first time a thorough evaluation and comparison of the mentioned methods at 
    hand of the complete TOSCA dataset.
\par
    {\em (iv)} We investigate here for the first time {\em soft correspondences} for use with the MOR signatures.
    The soft correspondences are derived from soft maps \cite{Eisenberger2018,SNBBG-2012}. At hand of dedicated
    evaluations we show that the MOR technique works very favourably with this technique to establish correspondences.
\par
    {\em (v)} Because of the significantly better stability that is achieved via the
    introduced numerical similarity transform, it turns out we gain {\em much more efficient} computational approaches
    for all the methods in virtually all of the considered settings, as the problem dimensions that are of interest can
    be reduced (except for the kernel-based methods when evaluating the geodesic error measure in some experiments where
    we do not observe much improvement).
\par    
Let us note that in addition, in all of these points we introduce novel aspects and techniques compared to \cite{BDB-2018}.
Especially, the extended experimental evaluations presented here allow to highlight some relevant properties 
of the MOR signature that have not been discussed earlier in \cite{BDB-2018}. As indicated in point (v) above,
our results also shed light on the number of eigenvalues that should be employed within all of the mentioned techniques
when using the numerical similarity transform as indicated above. This is an important practical issue for the discussed
methods for pointwise shape correspondence.

\paragraph{Paper Organisation}

After a broader exposition on related work on shape descriptors, we briefly recall the general
framework we rely on along with the arising PDEs and the numerical discretisations employed in space and time.
As our work relies very much on related numerical techniques, we give a comprehensive discussion of numerical solvers
with an emphasis on solvers for the occurring specific linear systems of equations.
In Section 6 we proceed with a detailed exposition of MOR methods and some related issues in our setting.
Our first experimental evaluation presented in Section 7 focuses on the evaluation of solvers in the context of MOR
techniques and time integration. This is followed by a discussion of the numerical similarity transform
and time rescaling we propose here for optimisation. Beginning with Section 9 we include
the kernel-based methods in our discussion. In Section 10
we give a detailed exposition of the comparison of the presented approaches, focusing thereby on the most promising
MOR methods identified in the previous sections, and including the soft correspondence techniques. The paper is finished
by a summary with a conclusion.

\section{More on Related Work}

A large variety of feature descriptors for 3D objects has been proposed in the fields of computer vision,
pattern recognition and computer graphics within the last decades.
The classic feature descriptors handle rigid transformations while more recent approaches on which we focus
are invariant under isometric transformations. Furthermore, one may identify in a broad sense four classes
as follows.

\paragraph{Distance-Based Methods}

In \cite{OFCD} a descriptor is proposed that uses the distribution of the Euclidean distances between pairs of randomly selected
points on the surface of a 3D model.
Other approaches use intrinsic properties of the shape by collecting information from intrinsic distances \cite{KS,SSKG}.
Later, a pose-invariant shape descriptor is introduced in \cite{GSC} as a 2D histogram estimated from 
the measure of the diameter of the 3D shape in the neighbourhood of each point on a surface and 
the average intrinsic distance from one point to all other points on the shape.
The authors in \cite{IAPKC} introduced a feature descriptor based on the distribution of the lengths of 
the longest geodesics on the shape. 
Finally, \cite{COO} introduced a multiscale signature based on the topological structure of the distribution of
geodesic distances centred at a given point.

\paragraph{Local Descriptors}

The idea of local descriptors is to characterise each point on a shape by using local geometric properties. 
Important classic examples of local point-based signatures are the well-known spin images \cite{JH} and shape context \cite{BMP}. 
Later, a multi-scale version of the local neighbourhood of the given point was used in \cite{YLHP} to construct a feature descriptor. 
The authors in \cite{PWHY} employed this approach for the computation of integral invariant features \cite{MHYS}. 
In \cite{ZBH}, the well-known SIFT descriptor \cite{Lowe-2004} which is based on locally computed features
is extended to 3D objects. A similar feature descriptor, called SHOT, which describes the local shape structure by a collection
of the distribution of normal directions, was introduced in \cite{TSS}.
Finally, among the newer methods, there is the Anisotropic Windowed Fourier Transform (AWFT) descriptor \cite{MRCB}. Starting from a 
collection of functions, the descriptors are obtained as a weighted linear combination of the coefficients of the AWFT.

\paragraph{Time-Evolution Methods}

Another trend in shape analysis consists of exploiting intrinsic time evolution processes carried by partial differential
equations on geometric shapes. In this context diffusion processes are well established, allowing a meaningful interpretation
relating the propagation of information and intrinsic distances.  
For example, the propagation of heat on a shape can be interpreted as a random walk among surface points \cite{CL,BB-2011}. 
In the spirit of this framework, \cite{SOG} introduced the HKS. The HKS describes the amount of heat
remaining at a certain point after a certain amount of time. The geometric interpretation
of this approach is that one can determine a connection between the heat kernel and intrinsic distances
via Varadhan’s formula \cite{Varadhan1967}. Later, a scale invariant extension of the HKS was developed \cite{BK-2010}.
In \cite{ASC} another feature descriptor namely the WKS inspired by equations of theoretical
physics is proposed. Based on the Schrödinger equation, the WKS represents the average probability of measuring a quantum
mechanical particle at a specific location. At the same time, \cite{BB-2011} proposed a scheme that is able to generalise
the diffusion-based approaches.

All of these methods are based on the spectral decomposition of the Laplace-Beltrami operator on manifold shapes. To this end, the
feature descriptors can be represented by a truncated series using the eigenvalues and eigenfunctions of the Laplace-Beltrami operator.
In \cite{MORCC} the authors developed a feature descriptor using a specific discrete diffusion process without truncating an
eigendecomposition.
As indicated in the introduction, some works have been proposed to avoid the need for the latter
by introducing an alternative method based on the
full discretisation of the underlying PDEs \cite{DBH-2017,DBH-2016} at the expense of additional computation time, when considering
all points on the shape.

\paragraph{Learning-Based Methods}

Feature descriptors based on intrinsic time-evolution methods are invariant to isometric transformations.
However, in applications elastic deformations may appear yielding some distortions in intrinsic distances.
In order to address this problem, a class of methods has been introduced by using the
learning-by-example approach. In \cite{LB-2014}, the so-called optimal point descriptor is proposed.
The technique proposed in \cite{KIDS} relies on random forests. 
In \cite{COC-2014}, a learning procedure has been introduced within the functional map framework,
and in \cite{LBBC} the authors proposed a learning strategy in the context oft the bag-of-word framework. 
Finally, \cite{BMMBCV} introduced a matching method that extends convolutional neural networks to the
spectral domain.

\paragraph{Some Remarks on Other Related Approaches}

Let us note that point or shape signatures can also be used within other important methods for shape correspondence that are extensions of classical point-to-point mappings.
As examples let us mention the works based on functional maps \cite{OBSBG} where
correspondences are modelled as linear operators between spaces of functions on manifolds.
In the technical context of functional maps let us also mention the use of product manifolds \cite{VLRBC}.
The functional map framework has been adopted and extended in several follow-up works \cite{ERGB,MRSWO,RCBTC,RLBBS}
with extensions to learning based techniques \cite{HLRBK,LRRBB} relating to the corresponding category as above.


\section{About the Shape Correspondence Framework}
In this section we introduce the basic facts that are necessary to define the shape correspondence framework.
Concerning the general shape analysis set-up, we largely follow concepts as discussed for
instance in \cite{TOSCA}.
For notions from differential geometry as employed here we refer the reader to \cite{DC}.
\subsection{Almost Isometric Shapes}
A three-dimensional geometric shape can be described by its bounding surface.
Thus, our shape model consists of compact two-dimensional Riemannian manifolds $\mathcal{M} \subset \mathbb{R}^3$, 
equipped with the metric tensor $g \in \mathbb{R}^{2\times 2}$ that describes locally the geometry.\par
Two shapes $\mathcal{M}$ and $\widetilde{\mathcal{M}}$ may be considered as
isometric if there is a smooth homeomorphism 
$T: \mathcal{M} \rightarrow \widetilde{\mathcal{M}}$ between the corresponding object surfaces
that preserves the intrinsic distances between surface points:
\begin{align}
d_{\mathcal{M}}(x_1,x_2) =  d_{\mathcal{\widetilde{M}}}(T(x_1),T(x_2)), \quad \forall x_1,x_2 \in \mathcal{M}
\end{align}
The intrinsic distance between two surface points $x_k$, $k=1,2$ may be interpreted as the shortest 
path along the surface $\mathcal{M}$ connecting $x_1$ and $x_2$.\par
In many applications, the notion of isometric shapes may be too restrictive. For instance,
small noise in a dataset could be considered as an elastic deformation yielding some distortions
in intrinsic distances. To take into account this issue, we call 
two shapes $\mathcal{M}$ and $\widetilde{\mathcal{M}}$ 
almost isometric, if there exists a transformation $S: \mathcal{M} \rightarrow \widetilde{\mathcal{M}}$
with
\begin{align}
 d_{\mathcal{M}}(x_1,x_2) \approx  d_{\mathcal{\widetilde{M}}}(S(x_1),S(x_2)),  \quad \forall x_1,x_2 \in \mathcal{M}
\end{align}
\subsection{PDE-based Models for Shape Description}
\label{PartialDifferentialEquationsonShapes}
A classic but still modern descriptor class that can handle almost isometric transformations relies on physical models
that are conveniently described by PDEs. In the following, we introduce the two fundamental PDEs that we employ to this end.
\par 
The {\em geometric heat equation} that yields a useful shape descriptor \cite{SOG} involves the Laplace-Beltrami
operator. This is the geometric version of the Laplace operator that takes into account the curvature
of a smooth manifold in 3D. Given a parametrisation of such a two-dimensional manifold,
the Laplace-Beltrami operator applied to a scalar function $u: \mathcal{M} \rightarrow \mathbb{R} $ can 
be expressed in local coordinates as
\begin{align}
 \Delta_\mathcal{M} \,  u = {\frac {1}{{\sqrt {|g|}}}} \sum_{i,j=1}^2 \partial _{i}\left({\sqrt {|g|}}g^{{ij}}\partial _{j} u\right)
\end{align}
where $g^{ij}$ are the entries of the inverse of the metric tensor and  $|g|$ is its determinant.
Using this the geometric heat equation reads as
\begin{align}
\partial_t u(x,t) =  \Delta_\mathcal{M} \,  u(x,t), 
\qquad 
x\in \mathcal{M}, \; t\in I
\label{heat} 
\end{align}
and describes how heat would diffuse along a surface $\mathcal{M}$.
\par
The {\em geometric wave equation} is the second PDE that is going to be discussed in this work.
It has been introduced in \cite{DBH-2017} as a useful model for computing 
pointwise a shape descriptor. Assuming the speed of wave propagation is identical to one in all
directions on the manifold, the corresponding PDE is
\begin{align}
\partial_{tt} u(x,t) =  \Delta_\mathcal{M} \,  u(x,t), 
\qquad 
x\in \mathcal{M}, \; t\in I
\label{wave} 
\end{align}
Both of the described PDEs require an initial condition in order to be meaningful. In the context of shape correspondence
construction, we employ a Dirac delta function $u(x,0)=u_0(x)=u_{x_i}$ centred around a point of interest $x_i\in\mathcal{M}$.
The PDE (\ref{wave}) is of second order in time, so that it needs to be supplemented not only 
by a spatial function as an initial state, but also an account of the
initial velocity of that initial state is needed. As it is a canonical choice, we consider the zero initial
velocity condition $\partial_{t}u(x,0)=0$. 
\par
Let us note that many shapes appear as a closed manifold with $\partial \mathcal{M}= \emptyset$, where 
it is not necessary to define additional boundary conditions. For the case $\mathcal{M}$ has boundaries, we may 
require $u$ to satisfy homogeneous Neumann boundary conditions.
\subsection{Feature Descriptor and Shape Correspondence}
We now make precise how geometric feature descriptors are obtained by employing the
introduced PDEs, and how we construct shape correspondence on that basis.

\paragraph{Feature Descriptor}
For many shape analysis tasks, it is useful to consider a pointwise feature descriptor 
as a shape representation. The purpose of the feature descriptor is to give an account 
of the geometry of the surface at a certain local region centred about a considered point.
To this end, we restrict the spatial component of solutions $u(x,t)$ of the introduced PDEs
to
\begin{align}
f_{x_i}(t):=u(x,t)_{\mid_{ x=x_i}} \quad \mbox{with}\quad  u(x,0)_{\mid_{ x=x_i}}=u_{x_i}\label{featureDescriptor}
\end{align}
and call the $f_{x_i}$ the feature descriptors at the location $x_i\in \mathcal{M}$.
\par
Let us comment that there exists a physical interpretation related to the feature descriptors
we employ. The heat based feature descriptor $f_{x_i}$ describes the rate of heat
transferred away from the considered point $x_i$. The spreading of heat takes into account the
geometry of the surface by using the Laplace-Beltrami operator.
In turn, the wave based feature descriptor
describes the motion amplitudes of an emitted wave front observed at the considered
point $x_i$ during time evolution. Therefore, the latter feature descriptor catches 
the typical wave interaction observable in the solution of the wave equation 
as can e.g.\ be seen in the well-known formula of d'Alembert in the 1D case.
Analogously to the situation for the heat-based descriptor,
the observable motion of the waves is influenced by the intrinsic geometry of the surface. 
As time evolves, the waves spread over the surface so that 
their amplitude observed via $f_{x_i}(t)$ have the tendency to decrease.\par
Let us note that the feature descriptors discussed in this work cannot distinguish 
between intrinsic symmetry groups as they rely on intrinsic shape properties.

\paragraph{Shape Correspondence}
To compare the feature descriptors for different locations $x_i\in \mathcal{M}$ and
$\widetilde{x}_j\in \widetilde{\mathcal{M}}$ on respective shapes 
$\mathcal{M}$ and $\widetilde{\mathcal{M}}$ , we
employ a distance $d_{f}(x_i,\widetilde{x}_j)$ using the $L_1$ norm as
\begin{align}
d_f(x_i,\widetilde{x}_j)= \int\limits_{I} | f_{x_i}-f_{\widetilde{x}_j}  | \, \mathrm{d}t
\end{align}
It is clear that the tuple of locations $(x_i,\widetilde{x}_j)\in\mathcal{M}\times\widetilde{\mathcal{M}}$ with the smallest
feature distance should belong together. This consideration naturally 
leads to a minimisation problem for all locations:
\begin{align}
(x_i,\widetilde{x}_j)=\argmin\limits_{\widetilde{x}_k \in \mathcal{\widetilde{M}}} \, d_f( x_i,\widetilde{x}_k) 
\end{align}
By using $\widetilde{x}_j=S(x_i)=x_i$, the map $S$ can pointwise be restored for all $x_i$. 
Let us comment that without further alignment it cannot be expected that the restored map $S$ is injective
or surjective, since the minimisation condition is not unique.

\section{Basic Discretisation of Continuous-scale Models}
\label{sec:3}
In this section we recall the basics of the discretisation of the PDEs we employ.
In order to prepare for later developments that are at the heart of the contributions of this
paper, let us note that we really give just the description of the fundamentals. The
concrete schemes used for computations are relying on the technical building blocks
we introduce here.

\subsection{Discretising in Space and Time}
The discrete surface representation for computations is given by a triangular mesh which we denote
as $\mathcal{M}_d = (P, T)$, cf.\ Figure \ref{mesh}.
The underlying point cloud $P := \{x_1 , \ldots, x_N \}$ contains a finite number of
vertices in terms of coordinate points.  The mesh is constructed by connecting the vertices $x_i$ so that one 
obtains triangular cells. The individual triangles $T$ contain the neighbourhood relations 
between corresponding vertices.
As visualised in Figure \ref{mesh}, we let $\Omega_i$ be the barycentric
cell volume surrounding the $i$-th vertex. 

Turning from space to time discretisation, we define time intervals $I_k=[t_{k-1},t_k]$ and
set $t_0=0$ for subdividing the complete integration time $[0,t_M]$.

\subsection{Finite Volumes: Semi-Discrete Form}
\label{FiniteVolumeSetUp}
Letting for the moment $\partial_*$ be either $\partial_t$ or $\partial_{tt}$, we consider the PDEs
(\ref{heat}) and (\ref{wave}) over a so-called control volume $\Omega_i$ and a time interval $I_k$. 
Integration in space and time yields
\begin{align}
\int \limits_{I_k}  \int \limits_{\Omega_i}    \partial_* u(x,t) \, \mathrm{d}x \,  \mathrm{d}t  
=  \int \limits_{I_k}  \int \limits_{\Omega_i}  \Delta_\mathcal{M} \,  u(x,t)   \, \mathrm{d}x \, \mathrm{d}t \label{he1}
\end{align}
In a finite volume method the quantities that are considered in computations are cell averages,
i.e.\ for the $i$-th control volume, or cell, 
we define
\begin{align}
u_i(t)=u(\bar{x}_i,t)= \dfrac{1}{|\Omega_i|} \int \limits_{\Omega_i}  u(x,t)  \,   \mathrm{d}x
\end{align}
where $|\Omega_i|$ denotes the area of the $i$-th control volume.
Therefore, we define the averaged Laplacian as
\begin{align}
L u_i(t)= \frac{1}{|\Omega_i|} \int \limits_{\Omega_i}   \Delta_\mathcal{M} \,u(x,t) \, \mathrm{d}x  
\end{align}
As for the meaning of the latter integral on the right hand side, one has to apply the divergence theorem
to substitute the volume integral into a line integral over the boundary of the cell volume.
\par
For discretisation of the arising integral quantities, we employ here the widespread cotangent weight
scheme as introduced in \cite{MDSB}, which we briefly recall now.

\begin{figure*}[t]
  \centering
  \begin{minipage}[t]{.47\linewidth}
    \centering
   \begin{tikzpicture}[]
\pgftext{ \includegraphics[width=0.9\textwidth]{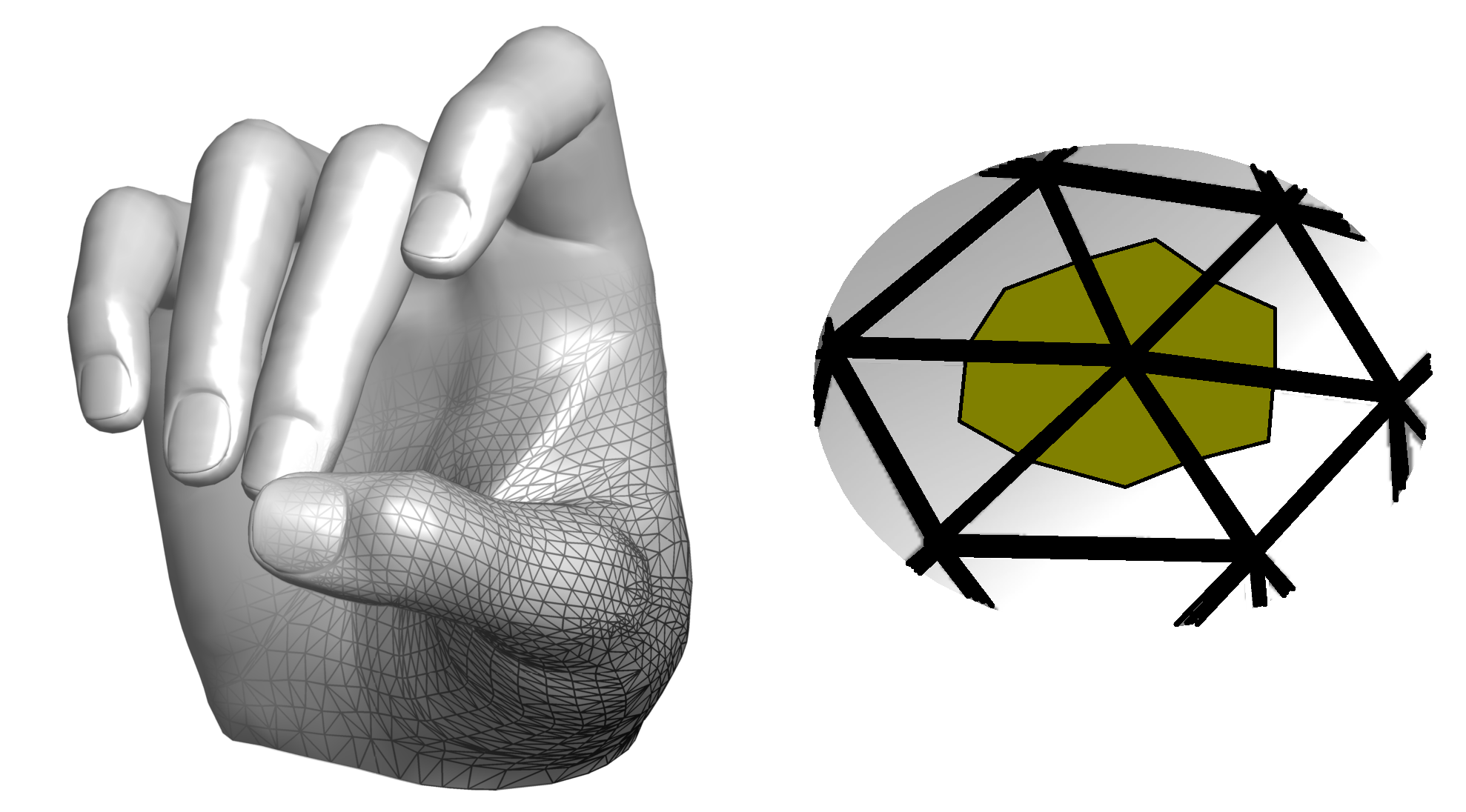} }  at (0pt,0pt); 
\draw [ultra thick] (1.97,0.1) ellipse (1.55cm and 1.24cm);
 \draw[ultra thick] (-1,0) -- (1,1.07);
\draw[ultra thick] (-1,-0.2) -- (1.3,-1.01);
  \end{tikzpicture}
\caption{Continuous and discrete shape representation (Figure adopted from \cite{BDB-2018}). The discrete shape is given by non-uniform linear triangles. Volume cells as shown here in green are constructed by using the barycentric area around a vertex.}
\label{mesh}
  \end{minipage}%
  \hfill%
  \begin{minipage}[t]{.47\linewidth}
    \centering
    \begin{tikzpicture}[]
\pgftext{ \includegraphics[width=1\textwidth]{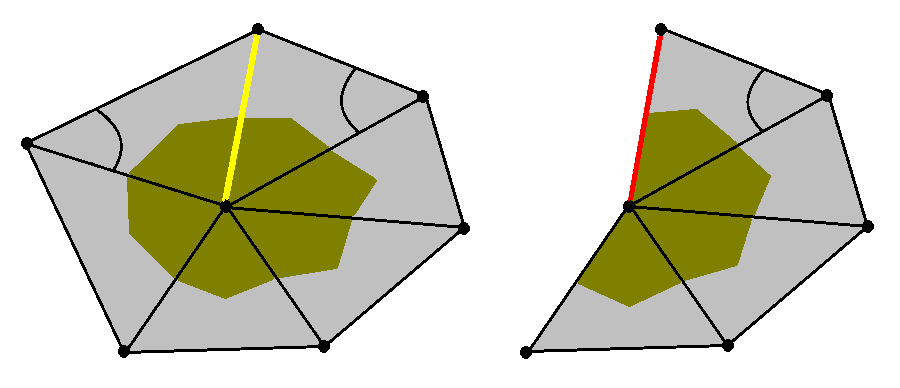} }  at (0pt,0pt);  
\node[text width=3cm] at (-2.0,0.5) {$\beta_{ij}$};
\node[text width=3cm] at (0.56, 0.85) {$\alpha_{ij}$};
\node[text width=3cm] at (4.25, 0.85) {$\alpha_{ij}$};
\node[text width=3cm] at (-0.62, -0.4) {$i$};
\node[text width=3cm] at (0, 1.6) {$j$};
\node[text width=3cm] at (3.05, -0.4) {$i$};
\node[text width=3cm] at (3.7, 1.6) {$j$};
\node[text width=3cm] at (-1.03, 0.5) {$w_{ij}$};
\node[text width=3cm] at (2.67, 0.5) {$w_{ij}$};
 \end{tikzpicture}
    \caption{The cotangent weight scheme as discretisation of the Laplace-Beltrami operator (Figure adopted from \cite{BDB-2018}).
Left: interior edge. Right: boundary edge.}
\label{weight}
\end{minipage}
\end{figure*}

The arising discrete Laplace-Beltrami 
operator $L\in \mathbb{R}^{N \times N}$ is composed of the sparse matrix representation
that can be written as $L=D^{-1} W$.
The appearing symmetric weight matrix $W$ contains the entries
\begin{align}
W_{ij}
= 
{\begin{cases} - \sum\limits_{j \in N_i}   
w_{ij},    
&{\mbox{if }} i=j  
\\  
\phantom{-}   w_{ij}, 
&{\mbox{if }} i \neq j  \text{ and } j\in N_i
\\
\phantom{-} 0, &{\mbox{else}}\end{cases}} 
\label{LBO}
\end{align}
where $N_i$ denotes the set of points adjacent to the vertex $x_i$.
The weights $w_{ij}$ of the edge $(i,j)$ between corresponding vertices distinguish between 
interior $E_i$ and boundary edges $E_b$ as shown in Figure \ref{weight}, and are given as
\begin{align}
w_{ij}
= 
{\begin{cases}
\displaystyle{\frac{\cot \alpha_{ij} + \cot \beta_{ij}  }{2 }     },
&{\mbox{if }} (i,j)\in E_i  
\\  
\displaystyle{\frac{\cot \alpha_{ij} }{2}     },
&{\mbox{if }} (i,j)\in E_b 
\end{cases}}
\end{align}
Furthermore, we let $\alpha_{ij}$ and $\beta_{ij}$ denote the two angles 
opposite to the edge $(i, j)$. The matrix 
\begin{align}
  D=\text{diag}\left(|\Omega_1|,\ldots,|\Omega_i|,\ldots,|\Omega_N| \right)
  \label{matrix-D}
\end{align}
contains the local volume cell areas.

Let us stress that $L$ is finally not symmetric. This fact has a significant influence on the computational
setting that will arise and which we discuss in detail later.

Lastly, we now put together and summarise the components of the discretisation as we developed it until now.
In this way we will end up with a semi-discrete form of the scheme.

Let a function defined on all cells be represented by now as an $N$-dimensional vector 
\begin{align}
  \mathbf{u}(t)=(u_1(t),\ldots ,u_N(t) )^\top
\end{align}
Rewriting \eqref{he1} using volume cell averages 
we obtain a system of ordinary differential equations (ODEs), one for each control volume:
\begin{equation}
\overset{*}{\mathbf{u}}(t) \, = \,    L \mathbf{u}(t) \label{ode}
\qquad
\textrm{where}
\qquad
\overset{*}{\mathbf{u}}(t)=\displaystyle{\frac{d^*\mathbf{u}(t) }{d^* t}}
\end{equation}
In order to clearly comment on the notation,
the use of the star derivative $\overset{*}{\mathbf{u}}$ indicates the time derivatives of first and second order,
respectively, as introduced in
\eqref{he1}.

Standard methods for the numerical solution of \eqref{ode} deal in general directly with first-order ODE systems.
For the geometric heat equation the system \eqref{ode} reads as
\begin{equation}
\mathbf{\dot{u}}(t) \, = \,    L \mathbf{u}(t) \label{odeHeat}
\end{equation}
In contrast, the geometric wave equation is of second order, and a transformation into a first-order system 
\begin{align}
  \mathbf{\dot{q}}(t) = K \mathbf{q}
  \label{odeWave}
\end{align}
with
\begin{align}
K= \begin{pmatrix}0 & I \\  L  & 0 \end{pmatrix} \in \mathbb{R}^{2N \times 2N}
\end{align}
where $I\in \mathbb{R}^{N \times N}$ is the identity matrix and 
\begin{align}
\mathbf{q}(t)=(\mathbf{q_1}(t),\mathbf{q_2}(t))^\top=(\mathbf{u}(t),\mathbf{\dot{u}}(t))^\top
\end{align}
is necessary.

Let us now turn to the discrete initial conditions of the described time evolutions.
The initial velocity function $\mathbf{\dot{u}}(0)$
for use with the geometric wave equation is identical to zero. Thus we
have only to describe here the discrete setting for the initial spatial
density $\mathbf{u}(0)$ which is used for both the geometric heat and wave equation.
To this end, we have to construct a discretised version of the Dirac delta function we formally
employed in the continuous-scale model.
Using the cell average 
\begin{align}
\int\limits_{\Omega_i}u(x,0) \, \mathrm{d}x  =1
\end{align}
where $u(x,0)$ may be interpreted now as a box function with unit area, being expressed as
\begin{align}
u(x,0)
= 
{\begin{cases}
\frac{1}{\abs{\Omega_i}},
&{\mbox{if }} x\in \Omega_i  
\\  
0,
&{\mbox{else }} 
\end{cases}}
\end{align}
the initial condition at the location $x_i$ can be formulated as 
\begin{equation}
\mathbf{u_{x_i}}=\mathbf{u}(x_i,0)=(0,\dots,0,|\Omega_i|^{-1},0,\dots, 0)^\top
\end{equation}
that implicitly bears a dependence on the index $i$. 
If we want to stress this dependence, we write the latter as $\mathbf{u_{i,0}}$.

%
%
%

\subsection{Discrete Time Integration}
\label{sec:13}
Solving the arising ODE systems \eqref{odeHeat} and \eqref{odeWave}, respectively, involves the application of numerical
integration and can be done by using time stepping methods. Widely used time integration schemes are the explicit Euler
method, the implicit Euler method and the trapezoidal rule known as Crank-Nicolson method.
\par
Explicit schemes require low computational effort when resolving one time step, but they are just conditionally stable and
may suffer from small time step restrictions, see e.g.\ \cite{Smith1985}. As we have shown in our conference paper
\cite{BDB-2018}, even recent explicit generalisations like the fast semi-iterative (FSI) method from \cite{HOWR2016} perform in an inefficient way
when applied for shape matching. The main reason is that in typical discrete meshes representing shapes, one has to face in general a large
variety of mesh widths, and especially also very small mesh widths arise. Since the spatial mesh width and the
allowed time step size of explicit methods are coupled, exactly this issue makes the explicit methods not attractive in our setting.
Let us note that this is the difference to the original setting in image processing from \cite{HOWR2016}, where the
typical image processing setting with a uniform mesh width of one is employed.
\par
In contrast, implicit methods like the implicit Euler method and the Crank-Nicolson method are
by theory unconditionally stable but they require to solve a system of linear equations.
Often the Crank-Nicolson method is used in applications due to its second-order convergence in time.
However, the latter method is not $L_0$-stable which may result in undesirable oscillations in the numerical solution
for problems with discontinuous initial conditions, compare \cite{Smith1985}. This happens to be the
setting of interest in our approach.
In such a case, usually $L_0$-stable schemes are preferred such that the implicit Euler method,
as \cite{DBH-2017,DBH-2016} has shown, represents a reasonable choice for our purpose.
\par
In the total, and as a consequence of the first investigations in \cite{BDB-2018}, we consider here only the implicit Euler method
for the numerical solution of the underlying geometric PDEs.
%


%
%
\paragraph{Implicit Euler Method for the Geometric Heat Equation}
First the application of the fundamental lemma of calculus for the left-hand-side of (\ref{odeHeat}) gives us
\begin{align}
\int \limits_{I_k}    \mathbf{\dot{u}}(t)  \,  \mathrm{d}t  =
\int \limits_{t_{k-1}}^{t_k}  \mathbf{\dot{u}}(t)  \,  \mathrm{d}t =  \mathbf{u}(t_k) -   \mathbf{u}(t_{k-1})
\end{align}
Subsequently, the approximation of the integral on the right-hand side of (\ref{odeHeat}) by using the right-hand rectangle method
\begin{align}
\int \limits_{t_{k-1}}^{t_k}  L \mathbf{u}(t)  \, \mathrm{d}t \approx \tau L \mathbf{u}(t_{k})
\end{align}
including the uniform time step $\tau=t_{k}-t_{k-1}$ and using the notation $\mathbf{u}(t_k)=\mathbf{u}^k$ finally results in
\begin{align}
(I- \tau L)\mathbf{u}^{k} =&  \mathbf{u}^{k-1} \label{BEheat}
\end{align}
with  $k\in \{1,\ldots, M\} $ and $\mathbf{u}^0=\mathbf{u_0}$. To compute the values $\mathbf{u}^{k}$ at time $k$ it requires
solving a large sparse system of linear equations in each time step.
\paragraph{Implicit Euler Method for the Geometric Wave Equation}
Analogous application of the same approximation scheme to the geometric wave equation leads to
\begin{align}
 \mathbf{q}^k = \mathbf{q}^{k-1} + \tau K\mathbf{q}^{k} \label{w0}
\end{align}
In \eqref{w0} the component $ \mathbf{q_1}$ at times $t_k$ and $t_{k-1}$ reads 
\begin{align}
\mathbf{u}^k &= \mathbf{u}^{k-1} + \tau \partial_t \mathbf{u}^{k} \label{w1}\\ 
\mathbf{u}^{k-1} &= \mathbf{u}^{k-2} + \tau \partial_t \mathbf{u}^{k-1} \label{w2}
\end{align}
while the component $ \mathbf{q_2}$ at $t_k$ can be written as
\begin{align}
 \partial_t \mathbf{u}^k = \partial_t \mathbf{u}^{k-1} + \tau L\mathbf{u}^{k} \label{w3}
\end{align}
The combination of \eqref{w1}-\eqref{w3} transform \eqref{w0} into a two step approach
\begin{align}
(I-\tau^2 L)\mathbf{u}^k = 2\mathbf{u}^{k-1} - \mathbf{u}^{k-2} \label{BEwave}
\end{align}
with a system size of $\mathbb{R}^{N \times N}$.

The geometric wave equation requires to define two initial conditions, namely $\mathbf{u}(t_0)$ and $\mathbf{\dot{u}}(t_{0})$.
With $\mathbf{u}^0=\mathbf{u_0}$ and the fixed initial velocity $\mathbf{\dot{u}}^0=0$
it follows
\begin{align}
  ( I - \tau^2 L)\mathbf{u}^1 =\mathbf{u}^0 \quad \textrm{for} \quad k=1
\end{align}

\section{Numerical Solvers for the Implicit Methods}

As described in the last section, the temporal integration will be done implicitly. In this case, there exist several
numerical solvers for the arising linear systems of equations, which have different advantages in terms of computational
effort and accuracy of the computed solution. In the following, we give a short overview of the methods that are useful
in our application.
\par
An essential key requirement for our objective of a correct shape matching is a sufficient accuracy of the computed
numerical solution. However, the underlying PDEs that are used to this end have to be solved for each point and on each
shape for a fixed time interval $t\in(0,t_M]$. Consequently, the computational costs are directly related to the number
of points of the regarded shapes. 
This suggests that one may forego high accuracy in exchange for a faster computational time. In order to evaluate this
proceeding, an analysis of the numerical solvers in context to shape matching is absolutely essential.
\par  
Let us start our discussion of implicit schemes by inspecting in more detail the arising linear systems.
The implicit schemes \eqref{BEheat} and \eqref{BEwave} result in a large sparse linear systems of equations and can
expressed as
\begin{equation}
A\mathbf{x}=\mathbf{b}\label{linS}  
\end{equation}
with $A=I- \tau L ,~\mathbf{b}=\mathbf{u}^{k-1},~\mathbf{x}=\mathbf{u}^{k}$ for the geometric heat equation and $A=I- \tau^2 L ,~\mathbf{b}=2\mathbf{u}^{k-1}-\mathbf{u}^{k-2},~\mathbf{x}=\mathbf{u}^{k}$ for the geometric wave equation.

Due to its nature as representing a discretisation of the
Laplace-Beltrami operator, the entries of $L$ show corresponding structure.
As a result the (for all time steps constant) matrix $A\in \mathbb{R}^{N \times N}$ is positive definite,
non-symmetric, large and sparse.

Implicit schemes are unconditionally stable without a time step restriction, however solving linear equations requires
significant computational effort and therefore a fast solver for large sparse linear systems of equations is necessary.
The linear system \eqref{linS} can be solved by using either sparse direct or sparse iterative solver.

\subsection{Sparse Direct Solver}
Application of sparse direct solvers, which are based on direct elimination as variations of the Gau\ss{} algorithm, are predestined for solving linear systems with a constant system matrix and multiple right-hand sides.
In that case, the underlying matrix $A$ will be factorised just once into a product of triangular matrices $A=LU$ by using a complete, sparse LU decomposition.
Subsequently, one can solve such systems for each right-hand side by forward and backward substitution, which is apparently highly efficient.

Direct methods are characterised by their extremely high accuracy of solutions. Their drawback is a high memory usage. Let us also note, that the computational costs will be at most $\mathcal O(N^2)$, where $N$ is the
number of equations.
\par
To improve the performance of the direct solver, an alternative object-oriented factorisation is useful to solve the linear system. In contrast to the LU decomposition, precomputing the matrix factors
in the object-oriented framework is more expensive, however this just has to be done once and it yields in the total a faster solver with exactly the same results. In order to accelerate the computation
in this way we use the SuiteSparse package \cite{Davis2013}.

\subsection{Sparse Iterative Solver}
In contrast, iterative methods are naturally not tweaked for extremely high accuracy, but they are very efficient in computing approximate solutions.
A particular class of iterative solvers that are used nowadays are Krylov subspace solvers, for a detailed exposition see e.g.\ \cite{Saad2003}.
The main idea behind the Krylov approach is to search for an approximate solution $\mathbf{x}_k$ of \eqref{linS} in a suitable low-dimensional
subspace $\mathbb{R}^{k}$ of $\mathbb{R}^{N}$, whereby the solution $\mathbf{x}_k$ is constructed iteratively.
The aim in the construction is thereby to have a good representation of the solution after a small number of iterations. Let us note that this
underlying idea is in general not directly visible in the formulation of a Krylov subspace method.

\paragraph{Conjugate Gradient Method}
The conjugate gradient (CG) method of Hestenes and Stiefel \cite{Hestenes1952} is probably the most famous Krylov subspace method and it is a
widely-used iterative solver for problems involving sparse symmetric and positive definite matrices. To realise the symmetric case
in \eqref{linS}, the multiplication from the left of the matrix $D$ in \eqref{matrix-D}
to the equations \eqref{BEheat} and \eqref{BEwave} results in
\begin{align}
(D- \tau W )\mathbf{u}^k &=D\mathbf{u}^{k-1}\label{BEheatSym} \\
(D- \tau^2 W)\mathbf{u}^k &=2D\mathbf{u}^{k-1}-D\mathbf{u}^{k-2}\label{BEwaveSym}
\end{align}
where $D-\tau^\alpha W$ with $\alpha=1,2$, is symmetric positive definite due to the
properties of the weight matrix $W$.

The CG method represents an orthogonal projection method and is based on the fact, that
solving the new system $\widetilde A\mathbf{x}=\mathbf{\widetilde b}$
with $\widetilde A=D- \tau W ,~\mathbf{\widetilde b}=D\mathbf{u}^{k-1}$ for the geometric heat equation and $\widetilde A=D- \tau^2 W ,~\mathbf{\widetilde b}=2D\mathbf{u}^{k-1}-D\mathbf{u}^{k-2}$ for the geometric wave equation,
can be reformulated as the minimisation of the quadratic function
\begin{equation} 
F(\mathbf{x})=
\frac{1}{2}\langle \mathbf{x}, \widetilde A\mathbf{x}\rangle_2-\langle \mathbf{\widetilde b},\mathbf{x}\rangle_2\label{CG-F}
\end{equation}
since
\begin{equation} 
\nabla F(\mathbf{x})=
\mathbf{0}
\quad
\Leftrightarrow
\quad
\widetilde A\mathbf{x}=\mathbf{\widetilde b}
\label{einschub-1}
\end{equation}
Thereby, $\langle \cdot,\cdot \rangle_2$ means the Euclidean 
scalar product.  It can be shown that the following property
\begin{equation}
\mathbf{x}_k\in \mathbf{x}_0+\mathcal{K}_k( \widetilde A,\mathbf{r}_0)\quad \mbox{with}\quad \mathbf{\widetilde b}-\widetilde A\mathbf{x}_k\perp \mathcal{K}_k(\widetilde A,\mathbf{r}_0)
\label{orthogonal}
\end{equation}
is fullfilled, whereby the $k$-th Krylov subspace $\mathcal{K}_k:= \mathcal{K}_k ( \widetilde A,\mathbf{r}_0)$ of $\mathbb{R}^N$ is defined as
\begin{equation} 
\mathcal{K}_k
:=
\mathsf{span}
(
\mathbf{r}_0,
\widetilde A\mathbf{r}_0,
\widetilde A^2\mathbf{r}_0,
\ldots
\widetilde A^{k-1}\mathbf{r}_0
)
\label{einschub-2}
\end{equation}
This means $\mathcal{K}_k$ is generated from an initial residual vector $\mathbf{r}_0= \mathbf{\widetilde b}-\widetilde A\mathbf{x}_0$ 
by successive multiplications with the system matrix $\widetilde A$. 

The condition \eqref{orthogonal} indicates that CG is an orthogonal projection method, and furthermore one can show that the
approximate solutions $\mathbf{x}_k$ are optimal in the sense that they minimise the so-called energy norm of the error vector. 
In other words, the CG method gives in the $k$-th iteration the best solution available in the generated subspace and the theoretical convergence is achieved at latest after the 
$N$-th step of the method \cite{Meurant1999}.

The main practical advantages of the CG method are the computationally cheap matrix vector multiplication in each iteration and the user-defined
termination if the approximate solution reaches a specific convergence tolerance. 
Typically, the \emph{relative residual}
\begin{align}
  \frac{\|\mathbf{\widetilde b}-\widetilde A\mathbf{x}\|_2}{\lVert\mathbf{\widetilde b}\rVert_2}\le \varepsilon
  \label{relative-residual}
\end{align}
is employed for defining the stopping criterion.
Increasing $\varepsilon$ leads naturally to faster computations but slightly worse results. The latter fact suggests to exchange high accuracy for a fast computational time.
Let us also note, that the total costs in the $k$-th iteration amount to at most $\mathcal O(kN)$.

One may also mention, that in practice one may suffer from convergence problems for very large systems, such that a \emph{preconditioner} is recommended to enforce all the beneficial
properties of the algorithm, cf. \cite{Benzi2002,Saad2003}. 
Moreover, it is usally required to fine-tune parameters in the preconditioned conjugate gradient (PCG) method.

\section{Model Order Reduction}
The introduced implicit methods have to handle large sparse systems, whereby the computational costs depends on the point cloud size.
{\em Model order reduction} (MOR) techniques can be used to approximate the full linear, time-invariant first-order ODE system \eqref{odeHeat}
and \eqref{odeWave} (where time invariance refers to the fact that the corresponding matrix is constant in time) by a very low dimensional
system, thereby preserving the main characteristics of the original ODE system. 
Existing MOR techniques can be classified in \emph{balancing based methods} and \emph{moment matching methods} respectively
in \emph{singular value decomposition based methods} and \emph{Krylov based methods}. The classification is not always consistent
in literature, for a general overview see \cite{Antoulas2005,Antoulas2015,Baur2014}.
\par  
In the following, we will describe the general procedure of MOR on the system \eqref{odeHeat} before discussing specific techniques,
obviously the approach is analogously applicable to \eqref{odeWave}.
\par
The mentioned methods are defined by a projected model which reduce the full system by projection matrices and rely on efficient
numerical linear algebra techniques for problems involving large sparse dynamic systems.
The basic concept of projection methods is to represent the high dimensional state space
vector $\mathbf{u}(t)\in\mathbb{R}^{N}$ in a reduced basis
\begin{equation}
\mathbf{u}(t)\approx V\mathbf{u_r}(t)\label{approximate}
\end{equation}
with $\mathbf{u_r}(t)\in\mathbb{R}^{r}$ and $r\ll N$, whereby $V\in\mathbb{R}^{N\times r}$ describes the projection matrix.
Applying this concept to \eqref{odeHeat} may be understood as projecting the original system 
\begin{equation}
     \begin{cases} 
     \begin{aligned}
     \mathbf{\dot{u}}(t)&=L\mathbf{u}(t) \\
         y_i(t)&=\mathbf{e_i}^\top\mathbf{u}(t), \quad \mathbf{u_i}(0)=\mathbf{u_{i,0}},~~i=1,\dots,N\label{ode2}
          \end{aligned}
     \end{cases} 
 \end{equation}
with the state variable $\mathbf{u}\in\mathbb{R}^{N}$, the single output variable $y_i(t)\in\mathbb{R}$ and the unit
vector $\mathbf{e_i}\in\mathbb{R}^{N}$ onto a reduced order model
 \begin{equation}
     \begin{cases}
     \begin{aligned}
     W^\top V \mathbf{\dot{u}_r}(t) &= W^\top LV\mathbf{u_r}(t) \\
        y_{r,i}(t)&=\mathbf{e_i}^\top V\mathbf{u_r}(t), \quad V\mathbf{u_{r,i}}(0)=\mathbf{u_{r,i,0}}\label{ode3}
        \end{aligned}
        \end{cases} 
 \end{equation}
by applying \eqref{approximate} and left multiplication with the projection matrix $W^\top\in\mathbb{R}^{r\times N}$.
The dynamical system in \eqref{ode2} can be interpreted as a zero-input-single-output (ZISO) system, whereby no input variable
exists in consequence of the considered boundary conditions. Moreover, we consider $N$ different initial conditions,
which corresponds to extracting the feature descriptor $y_i(t)$ for each point on the given shape.

By multiplication from left with $(W^\top V)^{-1}$, and assuming the inverse exists, the system $(\ref{ode3})$ leads to
 \begin{equation}
     \begin{cases} 
     \begin{aligned}
     \mathbf{\dot{u}_r}(t)& =(W^\top V)^{-1}W^\top LV\mathbf{u_r}(t)\\
      y_{r,i}(t)&=\mathbf{e_i}^\top V\mathbf{u_r}(t), \quad V\mathbf{u_{r,i}}(0)=\mathbf{u_{r,i,0}}\label{ode4}
        \end{aligned}
        \end{cases} 
 \end{equation}
 In general, the projection matrices are chosen as biortho\-normal matrices $W^\top V=I_r$ such that the reduced system of order $r$ can
 be described as follows:
 \begin{equation}
     \begin{cases}
     \begin{aligned}\mathbf{\dot{u}_r}(t) &=L_r\mathbf{u_r}(t) \\
        y_{r,i}(t)&=\mathbf{e_{r,i}}^\top \mathbf{u_r}(t), \quad V\mathbf{u_{r,i}}(0)=\mathbf{u_{r,i,0}}\label{ode5}
        \end{aligned}
        \end{cases} 
 \end{equation}
with $L_r=W^\top LV\in\mathbb{R}^{r\times r},~\mathbf{e_{r,i}}=\mathbf{e_i}^\top V\in\mathbb{R}^{r}$ and $\mathbf{u_{r,i,0}}\in\mathbb{R}^{r}$. \par
Let us mention, that projection techniques are characterised by the way of how to construct the projecting matrices $V$
and $W$. Common and widely used approaches in context of simulation of parabolic PDEs are modal coordinate reduction (MCR)
as employed in \cite{BDB-2018}, balanced truncation  (BT) method \cite{He2015}, proper orthogonal decomposition (POD)
\cite{Ojo2015} and Krylov subspace model order reduction (KSMOR) \cite{Sindler2013,He2015}. However, in our application
the usability of the BT and POD method can be excluded on grounds of efficiency in advance. Both methods are based on performing a
singular value decomposition, additionally BT has to solve Lyapunov equations and POD has to form the snapshot matrix for each point
on the given shape which results in inefficient processes.\par 
In the following, our aim is to give a short overview on the remaining abovementioned two methods, MCR and KSMOR.
However, as these approaches strongly rely on properties of the system matrix, we give a discussion
of the discrete Laplace-Beltrami operator beforehand.

\subsection{Discrete Laplace-Beltrami Operator}
\label{DiscreteLaplaceBeltramiOperator}
As indicated, we focus now on the properties of the discrete Laplace-Beltrami operator $L$.
More precisely, we consider the eigenvalues and eigenfunctions of $L$, which is equivalent to
considering the standard eigenvalue problem
\begin{equation}
 L\mathbf{v}=\lambda \mathbf{v}\label{EVP}
\end{equation}
whereby $L\in \mathbb{R}^{N \times N}$ is not symmetric. The last issue causes problems in terms of theoretical
and numerical aspects. On the one hand, non-symmetric matrices do not guarantee real eigenvalues and eigenvectors.
Secondly, their numerical computation may yield complex-valued results even if they were real. In the following,
we discuss these matters. A general overview on the Laplacian and its properties is presented in \cite{Zhang2010}.
\paragraph{The Generalised Eigenvalue Problem} As mentioned, the discrete Laplace-Beltrami operator can be written
as $L=D^{-1}W$, where $D$ is a regular diagonal matrix with positive entries on the diagonal, and $W$ is a symmetric
matrix. Under these conditions the eigenvalue problem \eqref{EVP} can be reformulated as a generalised eigenvalue
problem $D^{-1}W\mathbf{v}=\lambda \mathbf{v}$ or
\begin{equation}
 W\mathbf{v}=\lambda D\mathbf{v}\label{GEVP}
\end{equation}
which have the same eigenvalues and eigenvectors as the original problem. It should be noticed that if $W$ and $D$ are symmetric
and $D$ is also positive definite, which is the case here, then all eigenvalues $\lambda$ are real and the $N$ eigenvectors $\mathbf{v}$
are linearly independent, whereby the eigenvectors are $D$-orthogonal with $\mathbf{v_i}^\top D\mathbf{v_j}=\delta_{ij}$, see \cite{Parlett1998}.
This means, the eigenvectors are orthogonal with respect to the inner product 
\begin{equation}
 \langle f,g\rangle_D=f^\top D g\label{IP}
\end{equation}
Furthermore, the following equalities hold
\begin{equation}
 L=V\Lambda V^\top D,\quad I=V^\top DV,\quad \Lambda=V^\top WV\label{Properties}
\end{equation}
where $\Lambda$ is a diagonal matrix of eigenvalues and $V$ corresponds to the right eigenvector matrix of $L$.
Therefore, the eigenvalues of the underlying matrix $L$ are real and the eigenfunctions are $D$-orthogonal.\par
In addition, the weight matrix $W$ is a symmetric diagonally dominant matrix with real negative diagonal entries. From the Gershgorin circle
theorem it follows then that $W$ is negative semi-definite. The latter implies that $L$ is also negative semi-definite with respect to the inner
product \eqref{IP}, because of
\begin{equation}
 \langle \mathbf{x},L\mathbf{x}\rangle_D=\mathbf{x}^\top DL\mathbf{x}=\mathbf{x}^\top W\mathbf{x}\le 0
\end{equation}

\paragraph{A Note on Numerical Issues}
Although the eigenvalues of $L$ are real, the numerical solution of \eqref{EVP} may produce complex-valued results.
It is generally advantageous to compute $\mathbf{v}$ and $\lambda$ by making use of \eqref{GEVP} due to the fact, that numerical methods
for the generalised eigenvalue problem recognise the developed theoretical properties and consequently produce real eigenvalues and eigenvectors.

\subsection{Modal Coordinate Reduction (MCR)}
A suitable way to realise MOR in the context of the discussed application is the MCR
technique as used in \cite{BDB-2018}. This is a truncation approach which is based on the eigenvalue decomposition of the underlying
system matrix $L$ in \eqref{ode2}. The concept of MCR is to transform the full model from physical coordinates in physical space to modal
coordinates in modal space by using the eigenvectors of $L$ that are usually put together to form column by column an eigenvector matrix.
Subsequently, those modes are removed that have less important contributions to the system responses. Generally, only a few modes have a
significant effect on the system dynamics within the frequency range of interest.
\par 
The advantage of MCR is that the reduced model obtained in this way preserves the stability of the original system, however, the truncation may
not be optimal in the sense of the reduced eigenvalue spectrum of $L$. Let us briefly recall the approach.
\par
Considering again the semi-discrete scheme \eqref{ode2}, obviously a central issue is the Laplacian matrix $L\in\mathbb{R}^{N\times N}$.
Applying a regular (\emph{modal}) transformation $\mathbf{u}=V\mathbf{w}$ to the system, where $V\in\mathbb{R}^{N\times N}$ is the constant unit
eigenvector matrix of $L$, leads to
\begin{equation}
 V\mathbf{\dot{w}}(t) \, = \,    L V\mathbf{w}(t)
\end{equation}
Afterwards, the multiplications from the left by $D$ and $V^\top$ leads equivalently to
\begin{align}
V^\top D V\mathbf{\dot{w}} (t) =& V^\top DL V\mathbf{w}(t)\label{inv}
\end{align}
From \eqref{Properties} we have $I=V^\top DV$ and $\Lambda=V^\top WV$, whereby the latter is equivalent to $\Lambda=V^\top DLV$ due to $L=D^{-1}W$.
Inserting the last identities in \eqref{inv} results in
\begin{align}
\mathbf{\dot{w}} (t) =& \Lambda \mathbf{w}(t)
\end{align}
The latter equation is the starting point for choosing eigenvalues and eigenvectors (modes).
\par 
It is quite well-known, that the low frequencies which correspond to small eigenvalues are supposed to dominate the dynamics of the 
system. Suppose $r\ll N$ ordered eigenvalues $0=\abs{\lambda_1}<\abs{\lambda_2}\le\dots\le \abs{\lambda_r}$ 
are of interest. Consequently, we obtain with $\Lambda_r\in\mathbb{R}^{r\times r}$ extracted from $\Lambda$ and 
$V_r\in\mathbb{R}^{N\times r}$ the {\em reduced model} of order $r$
\begin{align}
  \mathbf{\dot{w}_r} (t) =& \Lambda_r \mathbf{w_r} (t)
\qquad
\textrm{where}
\qquad \mathbf{w_r}= V_r^\top D\mathbf{u}
  \label{odeRed}
\end{align}
This low dimensional system is much faster to solve than the original one. Applying the implicit Euler method to \eqref{odeRed} leads to  
\begin{align}
(I- \tau \Lambda_r)\mathbf{w_r}^{k} =&  \mathbf{w_r}^{k-1}, \quad\mathbf{w_r}^{0}=V_r^\top D\mathbf{u}^{0} \label{BEreduced}
\end{align}
or more precisely
\begin{align}
\mathbf{w_r}^{k} =&  P\mathbf{w_r}^{k-1}, \quad\mathbf{w_r}^{0}= V_r^\top D\mathbf{u}^{0} \label{BEreducedNew}
\end{align}
with
\begin{align}
  P=(I- \tau \Lambda_r)^{-1}=\text{diag}\Big(\frac{1}{1-\tau\lambda_1},\ldots,\frac{1}{1-\tau\lambda_r}\Big)
\end{align}
The reduced system \eqref{BEreducedNew} is based solely on diagonal matrices and can easily be solved by fast
sparse matrix-vector-multiplications.
\par
To summarise, the computational costs depend only on the eigenvalues and eigenvectors of $L$, which are
known to be computationally extremely intensive to obtain. However, on the basis of the above considerations only the smallest
eigenvalues have to be computed iteratively such that the MCR technique is practicable for a convenient number of modes.
\par
Let us mention, that MCR belongs to the projection based methods, because the reduced system in \eqref{BEreduced} can be gained
directly by using projection matrices $V$ and $W$ in \eqref{ode3} which correspond to the eigenvector matrix -- with an additional
dominant order -- of the Laplacian matrix $L$.


 \subsection{Krylov Subspace Model Order Reduction (KSMOR)}
 A popular tool in connection to MOR methods for dynamical systems are the Krylov subspace model order reduction (KSMOR) methods.
 These reduce  the original system in consideration of the input-output behaviour borrowing ideas from signal processing.
 In contrast to MCR, the KSMOR methods are based on  moment matching and approximate the transfer function of the original
 system \eqref{ode2}, which describes the dependence between the input and the output. The moment matching
 procedure relies on Krylov subspace methods and is therefore predestined for application in model order reduction
 of high-order systems. The main drawback of KSMOR methods is that in general there is no guarantee for
 preserving stability
 or error bounds of the original system. However, under certain assumptions it can be shown that the stability is preserved without
 additional effort, for instance this happens if the input matrix $L$ of the (stable) dynamical system is negative
 semi-definite \cite{selga2012}, which happens to be the case in this work. Let us briefly describe the approach.

 We now proceed by applying the Laplace transform at \eqref{ode2}, which transforms the original problem from time domain
 to frequency domain. Consequently, applying this transform to the ZISO system $\mathbf{\dot{u}}=L\mathbf{u}$ results in
 \begin{equation}
  sU(s)-\mathbf{u_{i,0}}=LU(s)\label{LT}
 \end{equation}
and further
 \begin{equation}
  U(s)=(sI-L)^{-1}\mathbf{u_{i,0}}\label{LT2}
 \end{equation}
 whereby the inverse $(sI-L)^{-1}$ exists for $s\neq \lambda_i$. The analogous procedure
 applied to the output $y_i(t)=\mathbf{e_i}^\top\mathbf{u}(t)$ leads to
 \begin{equation}
  Y_i(s)=\mathbf{e_i}^\top U(s)=\mathbf{e_i}^\top(sI-L)^{-1}\mathbf{u_{i,0}}\label{LT3}
 \end{equation}
 Contrary to the general representation, the output depends only on the initial condition.
 Subsequently, the transfer function is defined as
\begin{equation}
 g(s)=\mathbf{e_i}^\top(sI-L)^{-1}\mathbf{u_{i,0}}\label{transfer}
\end{equation}
which describes here the direct relation between the initial condition $\mathbf{u_{i,0}}$ and the output $Y_i(s)$ of the
original system in the frequency domain. The transfer function in \eqref{transfer} is called zero-input response and describes
the characteristics of the system itself. As already mentioned, the basic idea is to approximate the transfer function, whereby
the intended reduction focuses on matching the first moments. In particular, the moments $m_k(\sigma)=\mathbf{e_i}^\top (L-\sigma I)^{-(k+1)}\mathbf{u_{i,0}}$ are the coefficients of the Taylor series of the transfer function around $\sigma$
\begin{equation}
g(s)=-\sum_{k=0}^\infty m_k(s-\sigma)^k\label{taylor} 
\end{equation}
The latter is obtained in the following way: The reformulation of $g(s)$ in \eqref{transfer} leads to
\begin{align*}
 g(s)&=\mathbf{e_i}^\top(sI-L)^{-1}\mathbf{u_{i,0}}\\
 &=\mathbf{e_i}^\top\big((\sigma I-L)+(s-\sigma)(\sigma I-L)^{-1}(\sigma I-L)\big)^{-1}\mathbf{u_{i,0}}\\
 &=\mathbf{e_i}^\top\big(I+(s-\sigma)(\sigma I-L)^{-1}\big)^{-1}(\sigma I-L)^{-1}\mathbf{u_{i,0}}\\
 &=-\mathbf{e_i}^\top\big(I-(s-\sigma)(-\sigma I+L)^{-1}\big)^{-1}(-\sigma I+L)^{-1}\mathbf{u_{i,0}}
\end{align*}
and with denoting $\widetilde L=(-\sigma I+L)^{-1}$ we obtain
\begin{align*}
 g(s)=-\mathbf{e_i}^\top\big(I-(s-\sigma)\widetilde L\big)^{-1}\widetilde L\mathbf{u_{i,0}}
\end{align*}
which can be extend by using the Neumann series,
\begin{align*}
 g(s)&=-\mathbf{e_i}^\top\big(I+(s-\sigma)\widetilde L+(s-\sigma)^2\widetilde{L}^2+\ldots\big)\widetilde L\mathbf{u_{i,0}}\\
 &=-\sum_{k=0}^{\infty}\mathbf{e_i}^\top\widetilde {L}^{k+1}\mathbf{u_{i,0}}(s-\sigma)^k
\end{align*}
On this basis, the
aim is to determine a reduced system \eqref{ode5} such that its transfer function 
\begin{equation}
\widetilde{g}(s)=-\sum_{k=0}^\infty \widetilde{m}_k(s-\sigma)^k\label{taylor2} 
\end{equation}
 matches the first $q$ moments of the original transfer function i.e.
\begin{equation}
m_k(\sigma)=\widetilde{m}_k(\sigma), \quad k=0,\ldots,q-1\label{MM} 
\end{equation}
The moment matching technique can be done explicitly or implicitly. Explicit matching is known to be in general a numerically
unstable procedure, therefore implicit matching is performed based on Krylov subspaces.
\par
Proposed approaches to compute the Krylov subspaces of interest are variations of the Arnoldi and the Lanczos
process. In this work we apply the Arnoldi approach, which constructs an orthogonal basis
which we write as
\begin{equation}
\begin{aligned}
 V&=\mathcal{K}_q((L-\sigma I)^{-1},(L-\sigma I)^{-1}\mathbf{u_{i,0}})\quad \mbox{with}\\[4pt]
 \mathcal{K}_q&:=\mathsf{span} ((L-\sigma I)^{-1}\mathbf{u_{i,0}},\dots,(L-\sigma I)^{-q}\mathbf{u_{i,0}})\label{subspace}
\end{aligned}
\end{equation}
and where we especially employ $W=V$, such that $W^\top V=V^\top V=I_q$.
Let us note explicitly, that \eqref{subspace} means that the columns of $V$ are formed by the vectors constructed in $\mathcal{K}_q$.
It can be shown, that the first $q$ moments will
match if the projecting matrices $V$ and $W$ are bases of $\mathcal{K}_q$ in \eqref{subspace}.

We now choose $q=r$ from the exposition on MCR above, thus we choose the first $r$ moments for matching.
For $L_r=W^\top LV, ~\mathbf{e_{r,i}}=\mathbf{e_i}^\top V,~ W^\top \mathbf{u_i}=\mathbf{u_{r,i}}$ as well
as $(L-\sigma I)^{-1}\mathbf{u_{i,0}}=V\mathbf{r_0}$ with some $\mathbf{r_0}\in\mathbb{R}^{r}$, the first moment
of the reduced system reads as
\begin{align*}
&\widetilde{m}_0(\sigma)=\mathbf{e_{r,i}}(L_r-\sigma I_r)^{-1}\mathbf{u_{r,i,0}}\\
&=\mathbf{e_{i}}^\top V(W^\top LV-\sigma I_r)^{-1}W^\top \mathbf{u_{i,0}}\\
&=\mathbf{e_{i}}^\top V(W^\top LV-\sigma W^\top V)^{-1}W^\top \mathbf{u_{i,0}}\\
&=\mathbf{e_{i}}^\top V(W^\top LV-\sigma W^\top V)^{-1}W^\top (L-\sigma I)(L-\sigma I)^{-1}\mathbf{u_{i,0}}\\
&=\mathbf{e_{i}}^\top V(W^\top LV-\sigma W^\top V)^{-1}W^\top (L-\sigma I)V\mathbf{r_0}\\
&=\mathbf{e_{i}}^\top V(W^\top LV-\sigma W^\top V)^{-1}(W^\top LV-\sigma W^\top V)\mathbf{r_0}\\
&=\mathbf{e_{i}}^\top V\mathbf{r_0}=\mathbf{e_{i}}^\top (L-\sigma I)^{-1}\mathbf{u_{i,0}}=m_0(\sigma)
\end{align*}
The evidence of \eqref{MM} for the remaining moments can be done by induction, see \cite{Salimbahrami2002}.
Subsequently, the reduced system which has a small size, can be solved by using the direct solver and an incorporated
LU-decomposition. Let us mention, that the reduced system in \eqref{ode5} is in general not sparse.
\par
A parameter that still has to be determined is the choice of the expansion point $\sigma$. For
$\sigma=0,~ \sigma=\infty$ and $0<\sigma<\infty$ the resulting problem is known as Padé approximation, partial realisation and
rational interpolation, respectively. Let us recall at this point that after use of the Laplace transform we have to consider the
frequency domain. The value of $\sigma$ thus corresponds to the frequencies contained in the original model, such that
small values approximate low frequencies and $\sigma \rightarrow \infty$ higher frequencies. The underlying PDEs of this work
are characterised by a rather slow dynamic, therefore approximating the system at the frequency $\sigma\approx 0$
is a natural choice. In particular, for $\sigma=0$ the inverse $(L-\sigma I)^{-1}$ does not exist because $\lambda=0$ is an eigenvalue of $L$.
\par
For constructing the Krylov subspace $V=\mathcal{K}_q((L-\sigma I)^{-1},(L-\sigma I)^{-1}\mathbf{u_{i,0}})$, the included inverse
$(L-\sigma I)^{-1}$ leads again to the task of solving large sparse linear systems of equations. This requires the application of sparse direct
or sparse iterative solvers as introduced before. In this work, we apply the sparse direct solver in combination with the SuiteSparse
package.
\par
The computational costs of KSMOR are directly linked to the costs for construction of the Krylov subspace $V$.
However, at this point, it should also be mentioned once again that the projection matrix $V$ must be recalculated at each point and
on each shape due to the need to consider various initial conditions. 
Therefore, the computational costs scale substantially when increasing the number of Krylov subspaces.
Nevertheless, the KSMOR method can be highly practicable if a low number of subspaces represent a large part
of the system dynamics.

\section{Comparison of Solvers for Time Integration}
\label{Experiments:SolversforTimeIntegrationMethod}
As seen before, the temporal integration for \eqref{odeHeat} or \eqref{odeWave} will be done by using the implicit Euler method
which requires to solve systems of linear equations. For the latter there exist various numerical solvers, which have different properties
in terms of computational effort and accuracy of the computed solution which may clearly influence a shape matching based on it.
In the following, the numerical
solvers and their performance in terms of matching quality and run time will be analysed and evaluated for two different experiments.
The tests are only evaluated for the geometric heat equation \eqref{odeHeat}, analogous results have been achieved for \eqref{odeWave}
in undocumented tests. As methods of choice for solving the linear systems we consider:
\begin{enumerate}
\item The sparse direct solver. In addition, the SuiteSparse package is used.
\item The sparse iterative solver. For the termination of the CG method the parameter $\varepsilon>0$, which corresponds to the
  stopping criterion of the \emph{relative residual}, can be used as tuning parameter.
\item The MCR method. The performance of the solver can be tuned by the number of used eigenvalues and eigenvectors, here called modes $N_{\max}$.
\item The KSMOR method. The solver can be tuned by the number of used projection subspaces $\mathcal{K}_k$. For the computation of the
  subspace $V$, generated by the Arnoldi method, the sparse direct solver is used.
\end{enumerate}
\begin{figure}
  \begin{center}
   \includegraphics[width = 0.99\linewidth]{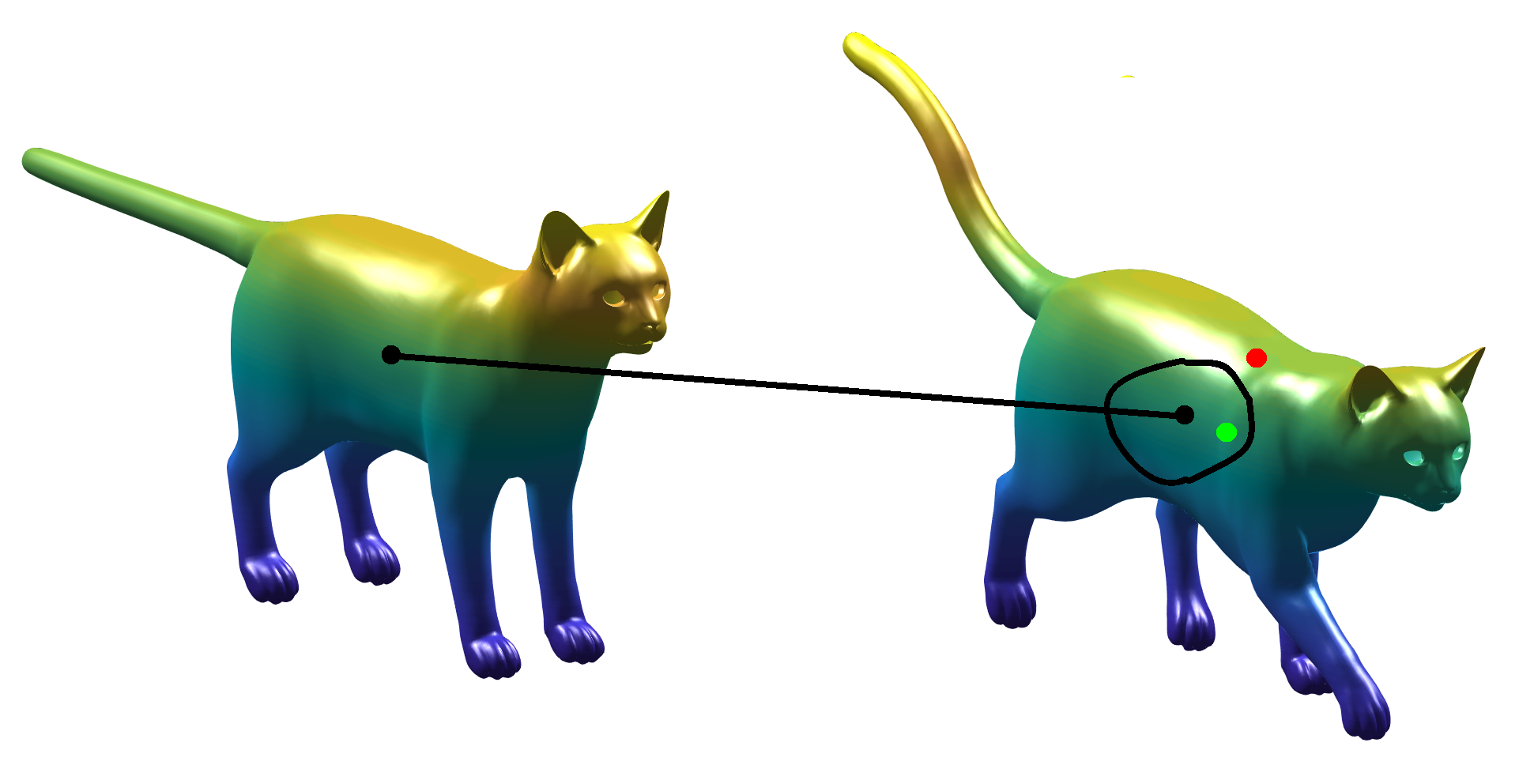}        
  \end{center}
  \caption{Evaluation of pointwise shape correspondence using the geodesic error. The (correct) ground truth matching is visualised
    by the black line. On the transformed shape we allow a certain radius of tolerance around the ground truth point. Matchings within
    the radius (e.g.\ the green point) are considered to be correct while points outside of the radius (e.g.\ the red point) are considered
    as outliers.}
      \label{geo}  
\end{figure}
As indicated, we are especially interested in run time as well as actual accuracy of the results in context of shape matching.
In order to evaluate the accuracy of the methods a dense point-to-point correspondence is performed, involving
all vertices the shapes are made off. In detail, the experiments are evaluated as follows.

\paragraph{Discrete Feature Descriptor}
The discrete version of the feature descriptor \eqref{featureDescriptor} is generated by numerical time integration of the underlying PDE. Therefore, the time axis $[0,t_M]$ is subdivided using $(M+1)$-time levels in $0:=t_0<t_1<\dots< t_M$  and the discrete feature descriptors have the form 
\begin{align}
\mathbf{f_{x_i}}=\big(u(x_i,t_0),\ldots,u(x_i,t_M)\big)^\top\in\mathbb{R}^{M+1}\label{discretefeatureDescriptor}
\end{align}
for $i=1,\dots, N$. The computation of a discrete feature descriptor at $x_i$ for a given shape requires to solve $M$ sparse linear systems of size ${N\times N}$. In total, the computation of all $\mathbf{f}_{x_i}$ implies that one has to solve in the total $N\cdot M$ linear systems for each shape.

\paragraph{Hit Rate} The percentage Hit Rate is defined as \\$TP/(TP+FP)$,
where TP and FP are the number of true positives and false positives, respectively.

\paragraph{Geodesic Error} For the evaluation of the correspondence quality,
we followed the Princeton benchmark protocol \cite{Kim2011}. This procedure evaluates the
precision of the computed matchings $x_i$ by determining how far are 
those away from the actual ground-truth correspondence $x^*$.
Therefore, a normalised intrinsic distance \\$d_\mathcal{M}(x_i, x^*)/ \sqrt{A_\mathcal{M}}$ on the transformed shape is introduced.
Finally, we accept a matching to be true if the normalised intrinsic distance is smaller than the threshold $0.25$, as illustrated in
Figure \ref{geo}.

\begin{figure}[h!]
  \begin{center}
    \begin{tabular}{cc}
      \includegraphics[width = 0.4\linewidth]{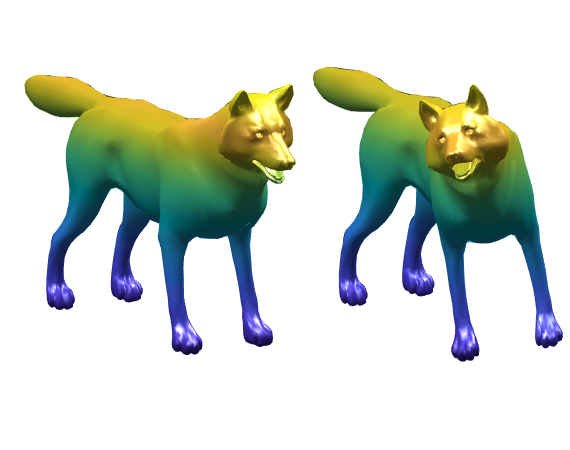}&\includegraphics[width = 0.4\linewidth]{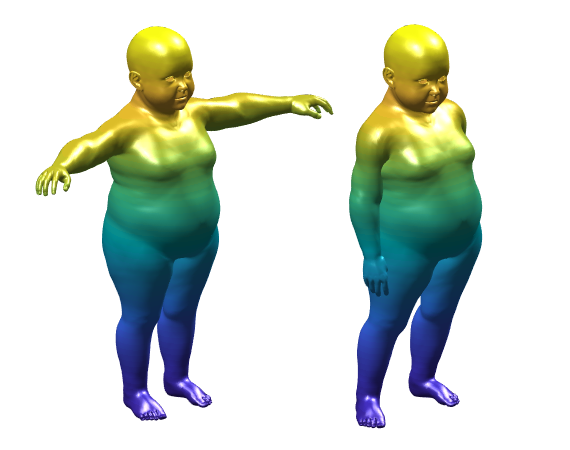}  \\
      Wolf & Baby \\
      $N=4344$ & $N=59727$
    \end{tabular}
  \end{center}
    \caption{For experimental evaluation, shapes at two different resolutions are compared, namely small and large. These are represented by the ``wolf" and ``baby" , taken from the
    TOSCA and KIDS dataset, respectively.}
      \label{dataset}  
\end{figure}

\paragraph{Dataset}
For the experimental evaluation, datasets at two different resolutions are compared, namely small and large. 
The small ($N=4344$) shapes of the wolf are taken from the TOSCA dataset \cite{TOSCA}. 
The baby shapes have a large resolution ($N=59727$) and are taken from the KIDS dataset \cite{KIDS}.
The datasets are available in the public domain, examples of it are shown in Figure \ref{dataset}.
All shapes provide ground-truth, and degenerated triangles were removed.

\paragraph{General Parameters}
We consider $M=25$ time levels such that the corresponding stopping time and the equidistant time increment are fixed to $t_M=25$ seconds
and $\tau=1$ second, respectively. The parameters are chosen without a fine-tuning, since we are interested to figure out the differences
of the numerical methods compared to accuracy and computational costs.
\par
Another issue by using time integration methods is the choice of the initial condition $u(x,0)=u_0(x)$.
As mentioned in Section \ref{FiniteVolumeSetUp}, we apply in this work a discrete delta peak in form of 
\begin{align}
\mathbf{u_{i,0}}=(0,\dots,0,|\Omega_i|^{-1},0,\dots 0)^\top\label{discreteIC}
\end{align}
\par
For constructing the Krylov subspace $V=\mathcal{K}_q((L-\sigma I)^{-1},(L-\sigma I)^{-1}\mathbf{u_{i,0}})$, the expansion point
is fixed to $\sigma=0.1$ without loss of generality.
\par
All experiments were done in MATLAB R2018b with an Intel Xeon(R) CPU E5-2609 v3 CPU. The computed eigenvalues
and eigenvectors for MCR are computed by the Matlab internal function \textit{eigs}.\par
Let us also note that the computations (in this experimental section) were taken by using the \emph{Parallel Computing Toolbox} 
integrated in Matlab. As mentioned before, all methods have to solve the geometric heat equation for each point of the 
given shape independently. This step can be easily parallelised by using the \emph{parfor} loop to distribute the code to 6 workers here. 

\subsection{Experimental Results}
\label{ExperimentalResults}

In what follows, we compare four important numerical approaches at hand of the two selected
test datasets.

\paragraph{Results on Sparse Direct Solver}
Firstly, we consider the wolf dataset with a small point cloud size of only $N=4344$ points for each of the two shapes, cf. Figure \ref{dataset}.
To compute the geodesic error, the geometric heat equation has to be solved numerically once for each point and on each shape. Ultimately, this
results in dealing with 217200 linear systems of size ${N\times N}$. By using the LU-decomposition offered by Matlab, the CPU time with around
600s offers significant running costs and is consequently quite inefficient. However, the computational costs can be reduced to around 25s by
using the powerful SuiteSparse package, which generates the same geodesic error accuracy.
\par
For the dataset baby, it is necessary to solve 59727 linear systems on each shape and for each time step. In this experiment the direct solver, including
SuiteSparse, performs this task in exactly 9267s ($\approx$ 2 1/2 hours). At this point it is recognisable, that large datasets produce high computational
costs and that the approach appears to be impractical for such shape matching applications.
\par
In following, the results (CPU time and geodesic error) of the direct solver with the object-oriented factorisation SuiteSparse are used for
the comparison of the remaining solvers.

\paragraph{Results on Sparse Iterative Solver}
\begin{figure}[b]
\centering
\includegraphics[width=0.48\textwidth]{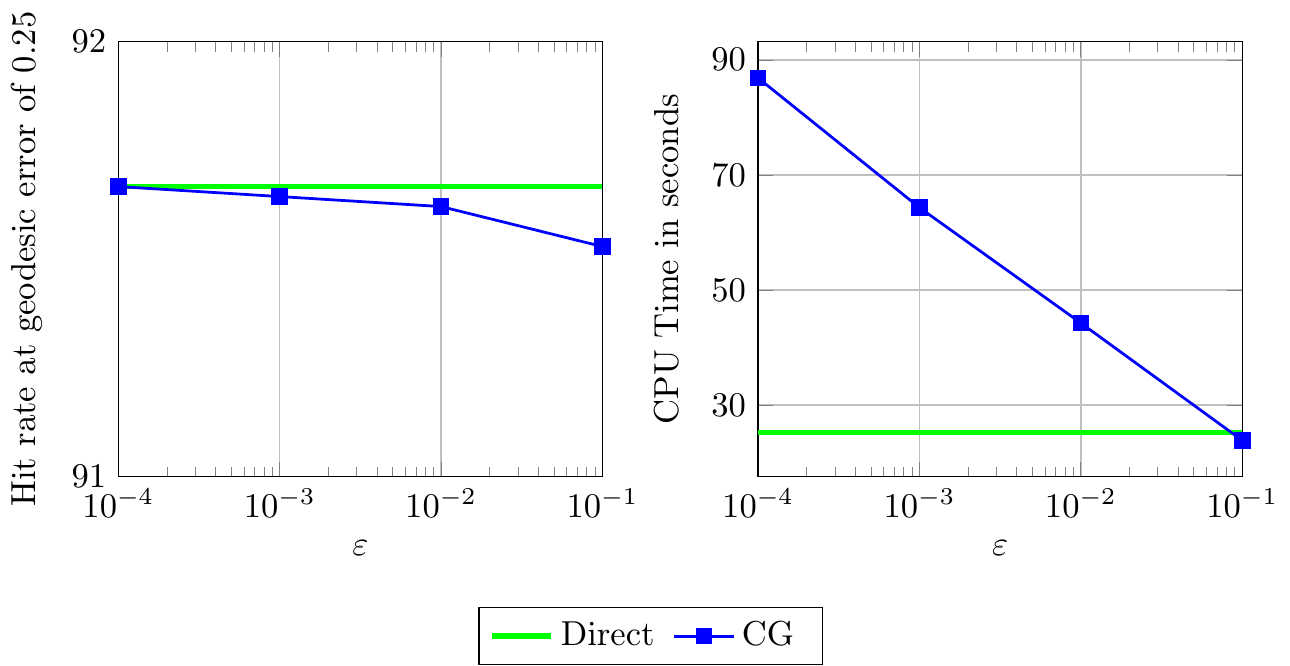}\\
\includegraphics[width=0.48\textwidth]{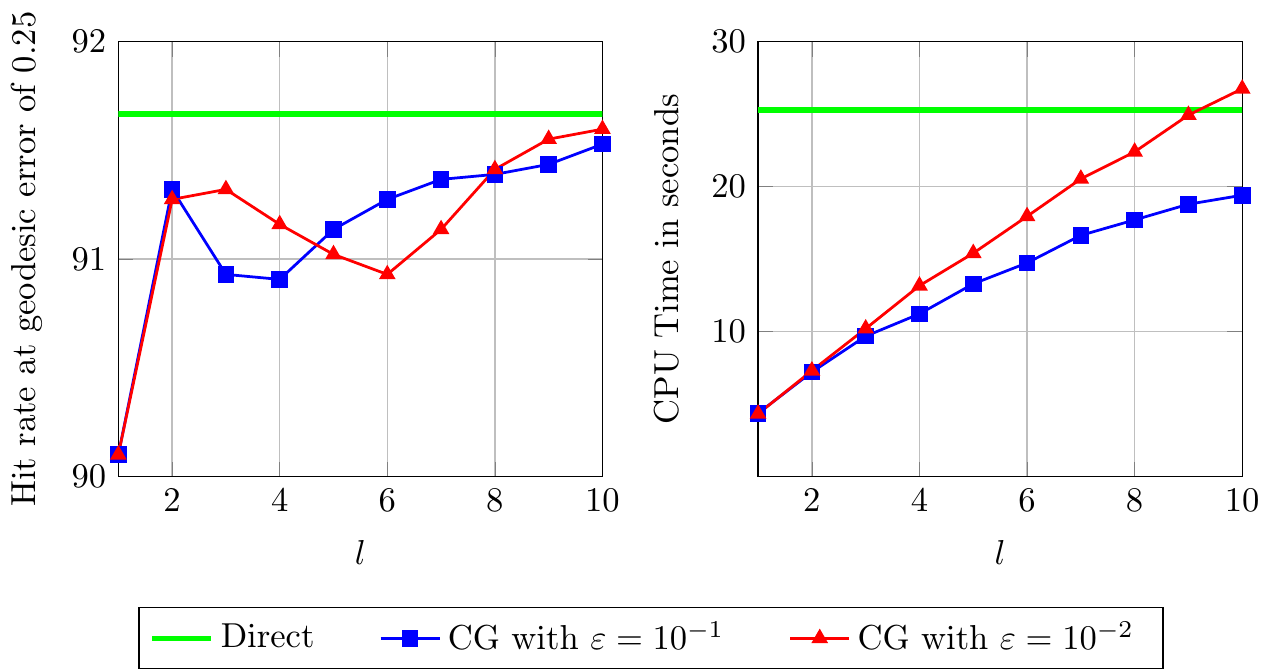}
\caption{Results on the dataset wolf by using the geometric heat equation. We compare the geodesic error at $0.25$ (left) and the performance time (right) between the direct solver and the CG method for various $\varepsilon$ (top) and different number of CG iterations 
$l$ (bottom) for $\varepsilon=10^{-1}$ and $\varepsilon=10^{-2}$.}
\label{wolf_CG}
\end{figure}
Solving the linear system \eqref{BEheatSym} by using the CG method involves fixing the stopping criterion of the iteration.
This involves setting the still free parameter $\varepsilon$ defining accuracy and one may define in addition an upper limit for the
number of iterations, compare also our first investigations in \cite{BDB-2018}.

Increasing $\varepsilon$ leads to a faster CG approximation, however it turns out that the accuracy remains almost unchanged also for the relatively large
value $\varepsilon=10^{-1}$, cf.\ Figure \ref{wolf_CG}. This means, a good representation of the solution is achieved already after a small number of
iterations. In addition, the repeated experiment for $\varepsilon=\{10^{-1}, 10^{-2}\}$ and various maximum numbers $l$ of CG iterations (i.e.\ the
iterative scheme terminates if one of the conditions is satisfied) is also shown in Figure \ref{wolf_CG}. 
In this case, the reduction of $l$ only has a minor effect on the geodesic error accuracy, yet it leads to a fast CPU time of only a few seconds.
\par
Regarding the baby dataset the matrix size should be taken into account. A typical problem of large systems is that the convergence becomes slower due
to the increase of the condition number of the system matrix, which is due to a larger system size. Therefore, we apply the PCG method including the
modified incomplete Cholesky (MIC) factorisation as a preconditioner, which is a modern variation of the standard incomplete Cholesky (IC) decomposition
and possesses in general a better performance. The mentioned factorisation uses a parameter-dependent, numerical fill-in strategy MIC($\gamma$),
where the parameter $\gamma$ is called \emph{drop tolerance} and describes a dropping criterion for matrix entries, cf.~\cite{Saad2003}.
In practice cf.\ \cite{Benzi2002,Baehr2017}, good results are obtained for values of $\gamma\in[10^{-4},10^{-2}]$, such that we fixed $\gamma=10^{-3}$.
\par
\begin{figure}[t]
\centering
\includegraphics[width=0.48\textwidth]{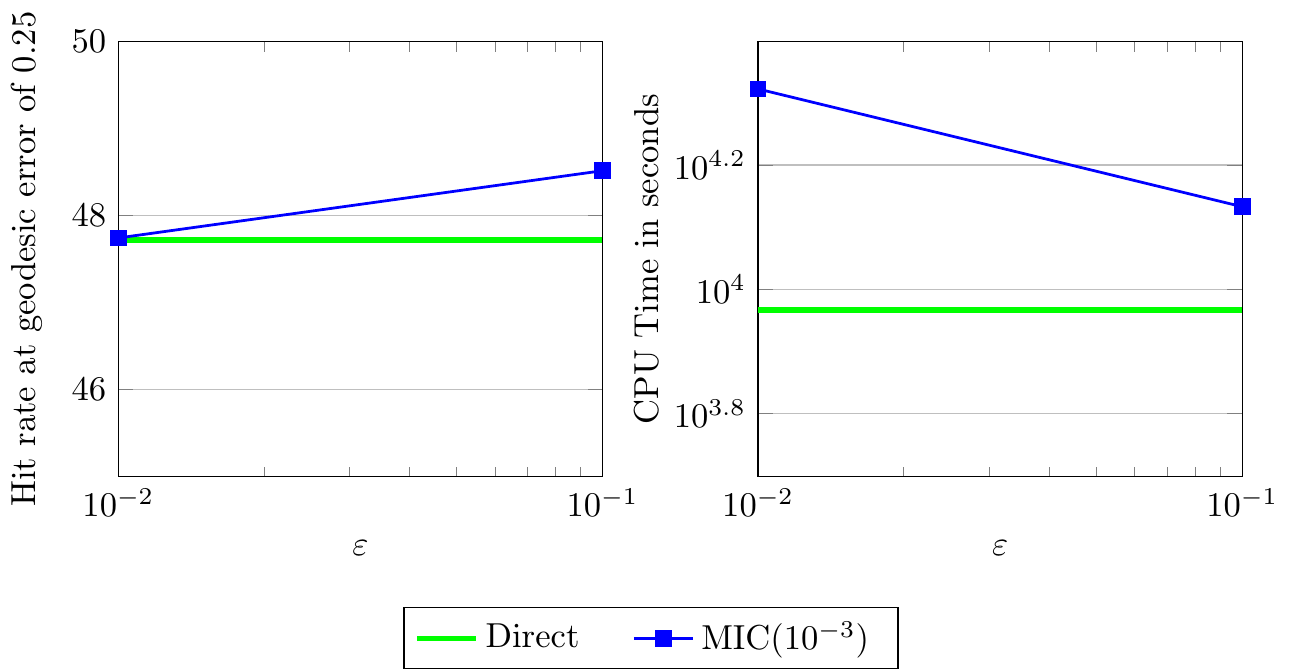}\\
\includegraphics[width=0.48\textwidth]{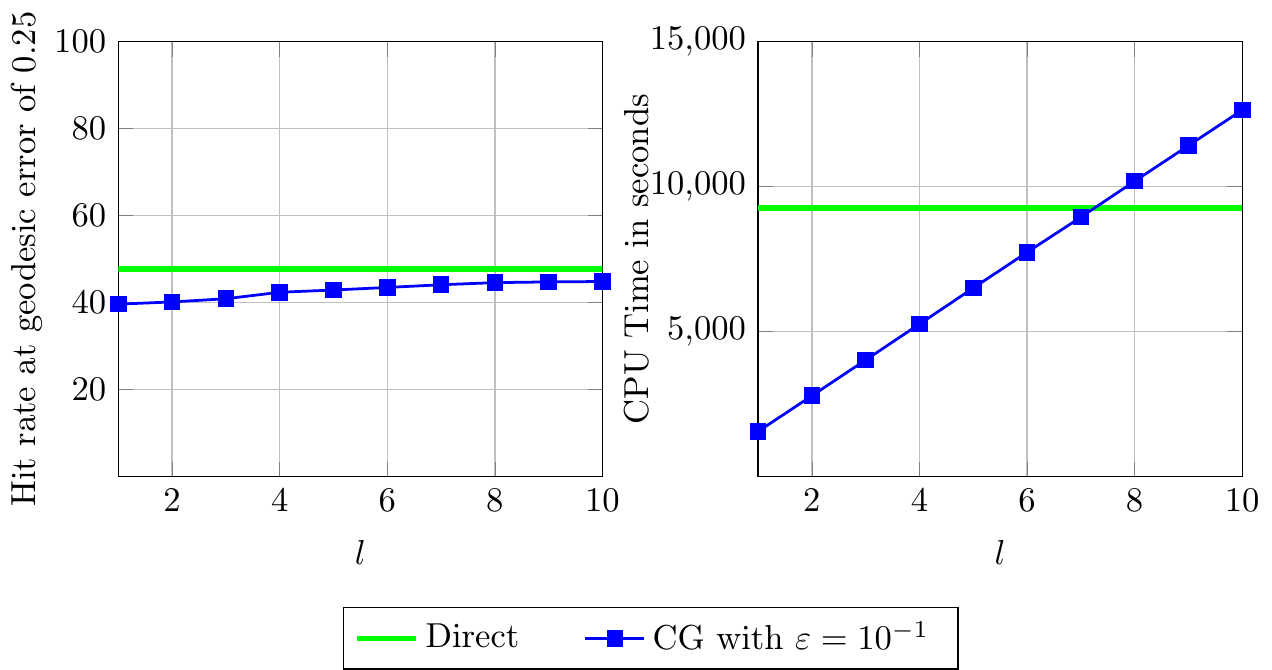}
\caption{Results on the dataset baby by using the geometric heat equation. We compare the geodesic error at $0.25$ (left) and the
  performance time (right) between the direct solver and the PCG method, including MIC$(10^{-3})$, for various $\varepsilon$ (top). In addition, we
  present a comparison between the direct method and the CG method (bottom) for $\varepsilon=10^{-1}$ and the first ten CG iterations $l$. The latter
  means, the iterative scheme terminates if one of the conditions is satisfied.}
\label{baby_CG}
\end{figure}

The results for MIC($10^{-3}$) are shown in Figure \ref{baby_CG}.
Compared to the performance of the direct solver no improvement is achieved. Furthermore, a closer examination of the required iterations for the
termination of the iterative algorithm shows that PCG needed just about 1 iteration. Therefore, PCG cannot be tuned further with regard to its performance.
The latter observation again inspires the idea to perform the CG method for a very small number of iterations $l\le 10$, which accordingly should be
sufficient to gain acceptable results in fast CPU time; see also our previous study \cite{BDB-2018} compared to which we present here some refinements
of the results. The results of CG for $\varepsilon=10^{-1}$ and small numbers $l$ of CG iterations are shown in Figure \ref{baby_CG}. Also in this case,
the iterative solver can reduce the computational effort significantly up to 80\%, however the percentage deviation of the accuracy in relation to the
direct solver is approximately up to 10\%. \par
\begin{figure}[t]
\centering
\includegraphics[width=0.48\textwidth]{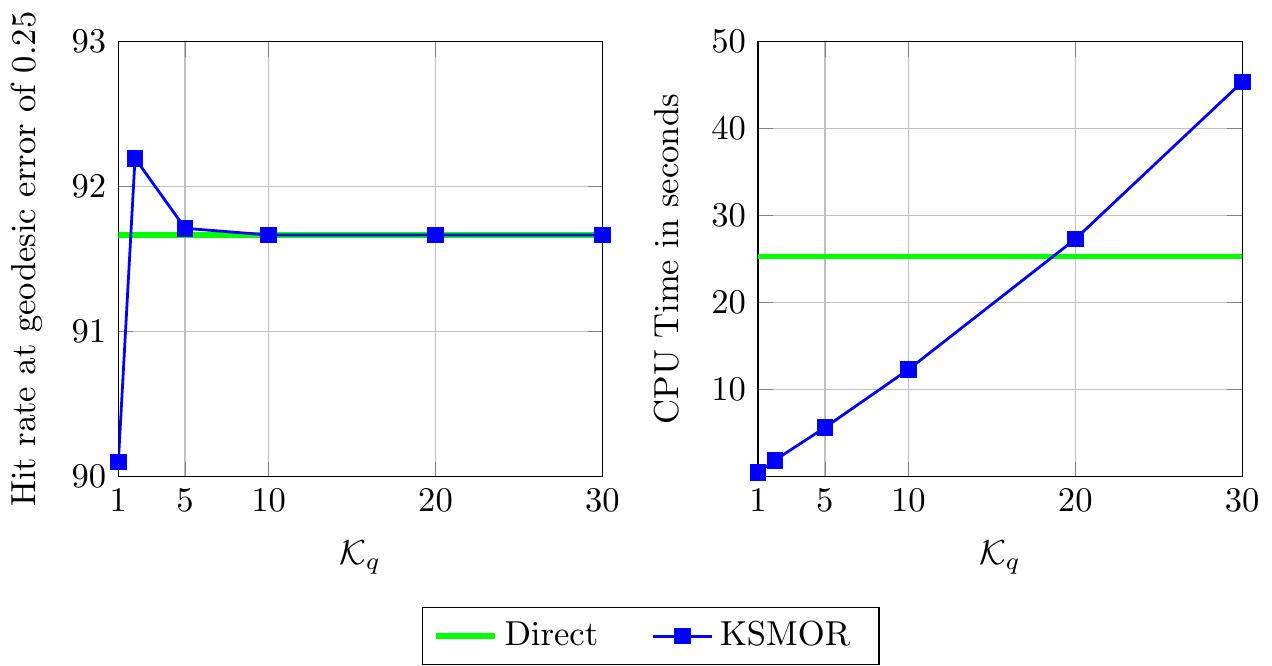}
\caption{Results on the dataset wolf by using the geometric heat equation. We compare the geodesic error at $0.25$ (left) and the performance
  time (right) between the direct solver and KSMOR method for a different number of Krylov subspaces $\mathcal{K}_q$.}
\label{wolf_Krylov}
\end{figure}

\paragraph{Krylov Subspace Model Order Reduction}
The results of KSMOR for the wolf dataset are illustrated in Figure \ref{wolf_Krylov}. Increasing the dimension of Krylov subspaces $q$ leads as expected, to a
more accurate approximation, whereby the same geodesic error accuracy compared to the direct solver can already be achieved for $q\le 10$. The required
computational costs correspond to those of CG, however we obtain a higher accuracy in terms of the geodesic error.
\par
\begin{figure}[t]
\centering
\includegraphics[width=0.48\textwidth]{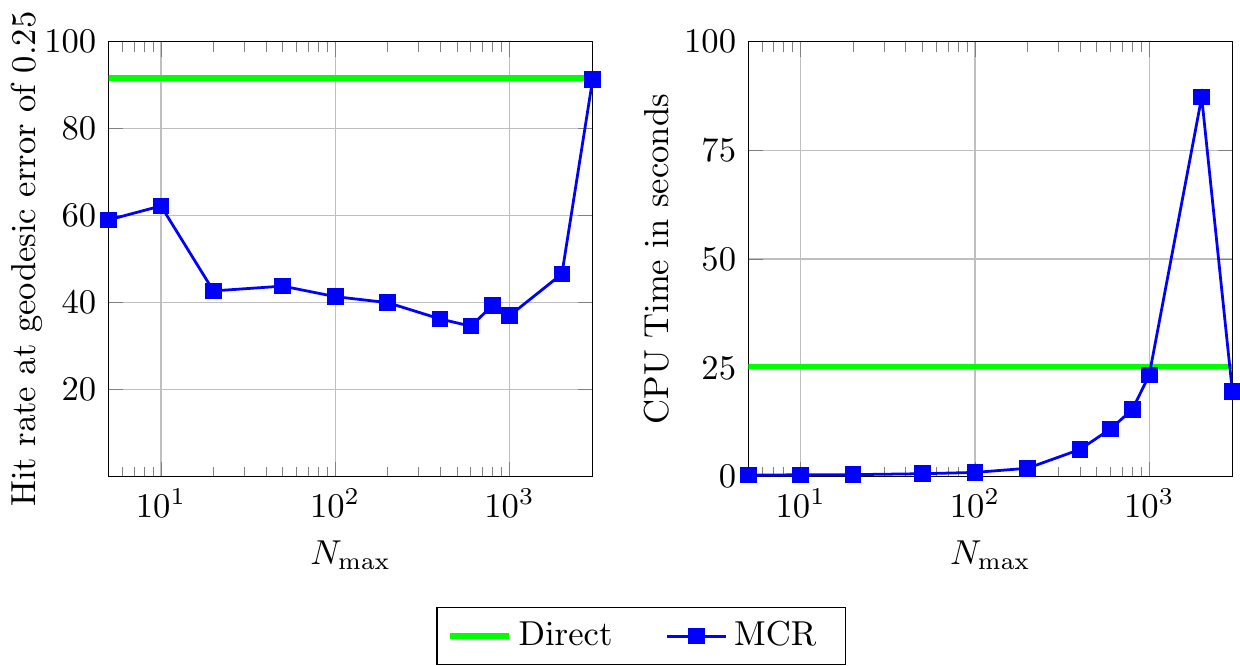}
\caption{Results on the dataset wolf by using the geometric heat equation. We compare the geodesic error at $0.25$ (left) and the
  performance time (right) between the direct solver and the MCR technique for different number of modes $N_{\max}\in[5,3000]$.}
\label{wolf_MCR}
\end{figure}

 
The same performance output is obtained by applying KSMOR for the baby dataset. Also for this dataset, using already approximately $q=2$ is deemed to be
sufficient to obtain an excellent trade-off between quality and efficacy and can save around 95\% of the computational time in relation to the direct solver.
However, the running time for $q=2$ with around 450 seconds is still significant for this standard size shape.
\begin{figure}[t]
\centering
\includegraphics[width=0.48\textwidth]{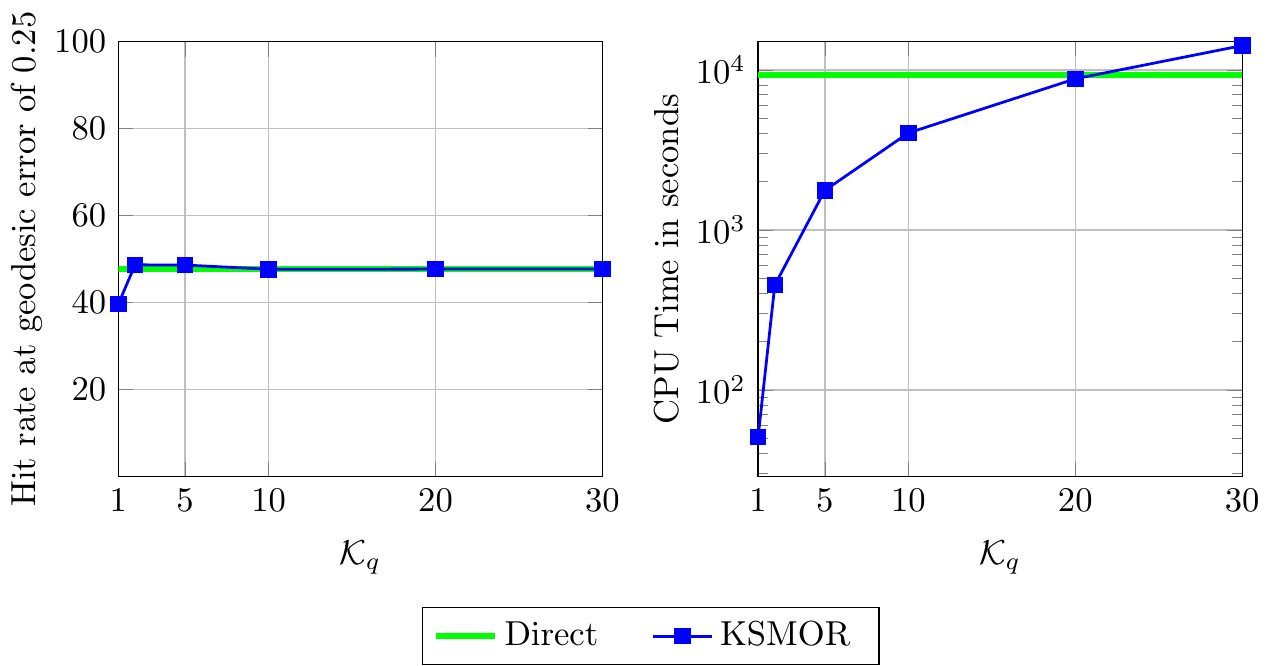}
\caption{Results on the dataset baby by using the geometric heat equation. We compare the geodesic error at $0.25$ (left) and the performance time
  (right) between the direct solver and KSMOR method for a different number of Krylov subspaces $\mathcal{K}_q$.}
\label{baby_Krylov}
\end{figure}
\paragraph{Modal Coordinate Reduction (MCR)}
Finally, we explore the MCR technique. The whole MCR process and reported computational times includes the computation of eigenvalues, eigenvectors
and solving of the subsequently reduced systems. For the experiments we increase the number of used ordered modes, starting from $N_{max}=5$ and going
up to $N_{max}=3000$.
Regarding to the wolf dataset we would expect that the correspondence quality gets better when increasing the number of modes. 
However, the evaluation in
Figure \ref{wolf_MCR} does not show this desirable behaviour. The outcome for a small spectrum $N_{\max}\approx 10$ are much better than for a
large spectrum $N_{\max}\approx 1000$, which seems at first glance not intuitive. The geodesic error oscillates in the range of $N_{max}\in[20,1000]$,
and beginning from about $N_{max}=2000$ which is a high number of modes for the approach
it converges against the solution of the direct method. Nevertheless, let us stress that the CPU time is incredibly fast when using a small number of modes.
For values up to ${N_{max}=100}$ the approximate solution is computed in less than 1 second, however the obtained geodesic error accuracy is not too high.
\par
\begin{figure}[t]
\centering
\includegraphics[width=0.48\textwidth]{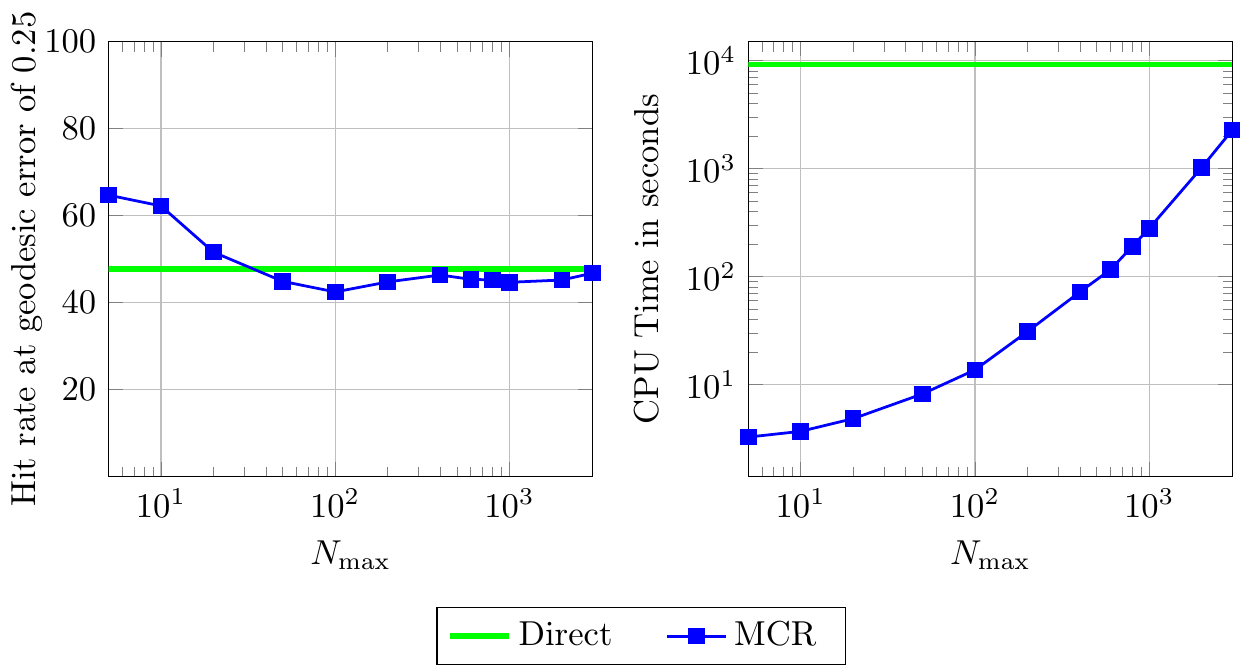}
\caption{Results on the dataset baby by using the geometric heat equation. We compare the geodesic error at $0.25$ (left) and the
  performance time (right) between the direct solver and the MCR technique for different number of modes $N_{\max}\in[5,3000]$.}
\label{baby_MCR}
\end{figure}

Applying the MCR technique to the baby dataset yields in the total a similar behaviour, but we also observe a surprising result, see Figure \ref{baby_MCR}.
In this experiment again, a small spectrum leads to better results compared to larger numbers of modes. In contrast to the wolf dataset, two
striking observations can be made. First, using $N_{max}=[5,10]$ even outperforms the geodesic error accuracy in relation to the direct solver and
secondly the geodesic error of both solvers are nearly equal already for $N_{max}=50$. 
In this case, MCR has an extremely short running time
of around 10 seconds, which reduces the computational effort compared to the direct solution by around 99.999\%.

\paragraph{Discussion of Solvers}
The experiments illustrated the performances of all applied solvers. The direct solver produces the best results in relation to the geodesic
error accuracy. However, for high resolution shapes ($N>50000$) the method is quite inefficient and generates computational costs amounting
to several hours.\par
The CG method may reduce the computational costs by around 80\%, whereby the percentage deviation of the accuracy in relation to the
direct method is approximately less than 10\%.
\par
Applying the KSMOR method achieves an even better performance. In this study, the choice of just $q=2$ saves around 95\% of the computational time,
whereby the percentage deviation of the geodesic error accuracy in relation to the direct solver is only around 1\%.
\par
The computational costs of the MCR method grow exponentially (by increasing $N_{max}$) and accordingly a practicable value $N_{max}$ should be small.
It is remarkable that the results for a small spectrum ($N_{max}\approx10$) are similar or even better to the ones for a significantly larger spectrum
($N_{max}\approx1000$). In case of a small spectrum, the MCR method is highly efficient and can save around 99\% CPU time, which suggests that this
technique is favourable for solving shape matching by time integration.
\par
Due to the incredible power of MCR, subsequently we follow this technique within this paper. Moreover, an interesting aspect would be to tune the
method so that it becomes more stable when using it with shapes such as the wolf shape, and in addition to improve the performance of the time
integration approach.
\par
Let us stress that we selected the two datasets used above since they are useful for demonstrating in detail the described behaviour of the schemes
as it can be observed in all tests we performed. They show the typical range of results that one obtains in terms of quality and efficiency.

\section{Optimised Numerical Signatures}
\label{OptimisedParameters}
As a result of the last section, the MCR method appears favourable for computing the numerical shape signature, especially when considering the extremely 
low computational costs. However, the geodesic error accuracy seems to depend on the given shapes and also on the number of used modes. For this
reason, the computation of the MCR signature leads to new challenges, as the aim is to ensure in a reliable way a high quality of results without
additional computational costs.

In the following, we introduce several means in order to enhance robustness and quality. In doing this, we significantly extend our previous
conference paper \cite{BDB-2018} as already indicated.

\subsection{Improved Eigenvalue Computation}
An important aspect of the MCR method is the way how to compute the eigenvalues $\lambda$ and eigenvectors $\mathbf{v}$ of the system matrix $L=D^{-1} W$.
As described in Section \ref{DiscreteLaplaceBeltramiOperator}, the computation is performed by solving the {\em generalised} eigenvalue problem \eqref{GEVP},
which uses certain properties of the underlying matrices $D$ and $W$ to produce numerically real-valued eigenvalues and eigenvectors.
However, the computation in this way may exhibit traces of instability in numerical practice.
\par
Another more favourable approach is to transform the generalised eigenvalue problem $W\mathbf{v}=\lambda D\mathbf{v}$ into the standard
eigenvalue problem without loss of the symmetry of $W$, such that more robust numerical eigensolvers can be applied.
The \emph{symmetric} standard eigenvalue problem can be achieved by a similarity transformation in the following way
\begin{align}
 L=D^{-1} W=D^{-\frac{1}{2}}D^{-\frac{1}{2}} WD^{-\frac{1}{2}}D^{\frac{1}{2}}=D^{-\frac{1}{2}}BD^{\frac{1}{2}}\label{similarityTransform}
\end{align}
whereby the matrix $B=D^{-\frac{1}{2}} WD^{-\frac{1}{2}}$ is symmetric due to the symmetry of $W$. Thus, the matrices $L$ and $B$ are similar
and have the same real eigenvalues. Moreover, the eigenpair $(\lambda,\mathbf{v})$ of $B\mathbf{v}=\lambda \mathbf{v}$ corresponds to the
eigenpair $(\lambda,\mathbf{\widetilde{v}})=(\lambda,D^{-\frac{1}{2}}\mathbf{v})$ of $L\mathbf{\widetilde{v}}=\lambda \mathbf{\widetilde{v}}$.
The matrix $D^{-\frac{1}{2}}$ in \eqref{similarityTransform} is well-defined due to the positive diagonal matrix $D^{-1}$ and is consequently also
positive diagonal. The eigenvectors $\mathbf{v}$ of $B$ are orthonormal with $\mathbf{v_i}^\top \mathbf{v_j}=\delta_{i,j}$ so that the
eigenvectors of $L$ are $D$-orthogonal
\begin{align}
\mathbf{v_i}^\top \mathbf{v_j}=(D^{\frac{1}{2}}\mathbf{\widetilde{v}_i})^\top D^{\frac{1}{2}}\mathbf{\widetilde{v}_j}=\mathbf{\widetilde{v}_i}^\top D\mathbf{\widetilde{v}_j}=\delta_{i,j}\label{IP2}
\end{align}
Let us briefly comment on an important issue here that may arise when adapting our framework to different types
of spatial discretisations such as e.g.\ the quite prominent finite element method which is one of the frequently used
methods for numerical computations. The proposed transformation is applicable whenever $D^{-1}$ (or $D$) is symmetric positive definite.
This situation will for instance occur when using the finite element ansatz instead of the presented finite volume set-up in
Section \ref{FiniteVolumeSetUp}, whereby $D^{-1}$ consequently represents the \emph{mass} matrix in that setting cf.\ \cite{Zhang2010}.
In this case, the generalised eigenvalue problem $W\mathbf{v}=\lambda D\mathbf{v}$ can be reformulated by using the existing
Cholesky decomposition $D=LL^\top$, with a lower triangular matrix $L$ and positive diagonal entries, as
\begin{align}
&W\mathbf{v}=\lambda D\mathbf{v}\nonumber\\
\Leftrightarrow\quad  &W\mathbf{v}=\lambda LL^\top\mathbf{v}\nonumber\\
\Leftrightarrow\quad   &L^{-1}W\mathbf{v}=\lambda L^\top\mathbf{v}\nonumber\\
\Leftrightarrow\quad   &L^{-1}WL^{-\top}L^{\top}\mathbf{v}=\lambda L^\top\mathbf{v}\nonumber\\
\Leftrightarrow\quad   &L^{-1}WL^{-\top}\mathbf{\widetilde{v}}=\lambda \mathbf{\widetilde{v}}\nonumber\\
\Leftrightarrow\quad  & B\mathbf{\widetilde{v}}=\lambda \mathbf{\widetilde{v}}
\label{choleskyTransform}
\end{align}
with $B=L^{-1}WL^{-\top}$ and $\mathbf{\widetilde{v}}=L^\top\mathbf{v}$. The eigenvectors are then retransformed according
to $\mathbf{v}=L^{-\top}\mathbf{\widetilde{v}}$. Clearly, for a diagonal matrix as in our method holds $D=D^{\frac{1}{2}}D^{\frac{1}{2}}$
so that \eqref{choleskyTransform} corresponds to \eqref{similarityTransform}.
\par
\begin{figure}[t]
\centering
\includegraphics[width=0.48\textwidth]{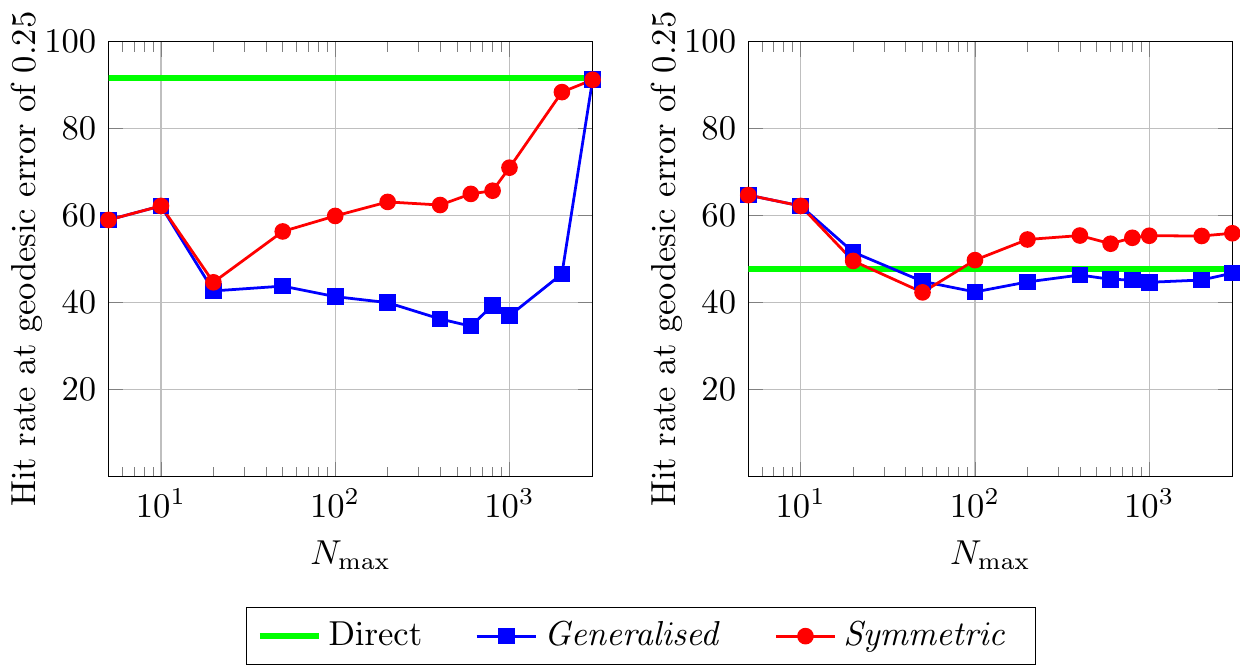}
\caption{Results on the dataset wolf (left) and baby (right) by using the geometric heat equation. We compare the geodesic error at $0.25$ between
  the direct solver and the MCR technique for different number of modes $N_{\max}\in[5,3000]$. The eigenvalues and eigenvectors are computed by
  solving the {\em generalised} and the \emph{symmetric} eigenvalue problem. Using the similarity transformation achieves more stable matching results.}
\label{gen_sim_MCR}
\end{figure}

In the following, we analyse the abovementioned computation of the eigenpairs $(\lambda,\mathbf{v})$ using the wolf and baby dataset. The evaluation
by solving the generalised and the symmetric eigenvalue problem is presented in Figure \ref{gen_sim_MCR}. As we indicated above, we observe that the
ansatz of using the similarity transformation for tackling the eigenvalue problem achieves significantly more stable matching results.
Let us emphasise that we see in the case of the wolf shape after this improvement the accuracy behaviour that we would have expected when
increasing the number of modes which indicates the better robustness of the approach, and also for the baby shape we obtain a quality improvement.

\subsection{Improved Scaling of the Integration Domain $[0,\mathit{t_M}]$}
We also discuss now the impact of the choice of the parameter $t_M$ with respect to the MCR technique.
The computed numerical signatures are computed over a temporal domain $(0,t_M]$, which makes them easily measurable for the shape matching purposes.
\par
Let us start our discussion of time scaling by considering some ideas for motivation. 
It is well-known, that the analytic solution of the semi-discrete geometric heat equation $\mathbf{\dot{u}}=L\mathbf{u}$ has the form 
\begin{align}
\mathbf{u}(t)=e^{Lt}\mathbf{u}(0)=e^{V\Lambda V^\top Dt}\mathbf{u}(0)=Ve^{\Lambda t}V^\top D\mathbf{u}(0)\label{solutionHeatMatrix}
\end{align}
or more precisely
\begin{align}
\mathbf{u}(t)=\sum_{i=1}^N e^{\lambda_i t}\mathbf{v_i} \mathbf{v_i}^\top D\mathbf{u}(0)=\sum_{i=1}^N\mathbf{v_i}^\top D\mathbf{u}(0) e^{\lambda_i t}\mathbf{v_i}\label{solutionHeatVector}
\end{align}
where $\mathbf{v_i}^\top D\mathbf{u}(0)$ is a scalar. Therefore, the pointwise feature descriptor has the form
\begin{align}
f_{x_j}(t):=u(x,t)_{\mid_{ x=x_j}}&=\Big(\sum_{i=1}^N\mathbf{v_i}^\top D\mathbf{u_j}(0) e^{\lambda_i t}\mathbf{v_i}\Big)^\top\mathbf{\hat{e}_j}\nonumber\\
&=\Big(\sum_{i=1}^N\mathbf{v_i}^\top \mathbf{\hat{e}_j} e^{\lambda_i t}\mathbf{v_i}\Big)^\top\mathbf{\hat{e}_j}
\label{solutionHeatVector2}
\end{align}
with 
$D\mathbf{u_j}(0)=\mathbf{\hat{e}_j}$, where $\mathbf{\hat{e}_j}$ is the $j$-th unit vector.
The representation of the solution \eqref{solutionHeatVector2} describes an exponential decay of heat transferred away from the considered point $x_j$.
Using only a small number of $r$ dominant frequencies (small eigenvalues) of a system, the reduced basis solution as approximative feature descriptor is given by
\begin{align}
f_{x_j}(t)\approx \widetilde{u}(x_j,t)=\Big(\sum_{i=1}^r\mathbf{v_i}^\top \mathbf{\hat{e}_j} e^{\lambda_i t}\mathbf{v_i}\Big)^\top\mathbf{\hat{e}_j}\label{solutionHeatVector3}
\end{align}
Furthermore, the reduced solution \eqref{solutionHeatVector3} with eigenvalues $0<\abs{\lambda_i}\ll 1$ ($i=1,\ldots, r$), can be simplified for small
times $t$ as
\begin{align}
\widetilde{u}(x_j,t)&=\Big(\sum_{i=1}^r\mathbf{v_i}^\top \mathbf{\hat{e}_j} \underbrace{e^{\lambda_i t}}_{\approx 1}\mathbf{v_i}\Big)^\top\mathbf{\hat{e}_j}\nonumber\\&\approx (V{\mathbf{y_j}})^\top\mathbf{\hat{e}_j}=c_j\label{solutionHeatVector4}
\end{align}
for all time levels $t=t_k$, where $c_j$ is a constant and ${\mathbf{y_j}}$ is the $j$-th column vector of $V^\top$.
\par
\begin{figure}[t]
\centering
\includegraphics[width=0.38\textwidth]{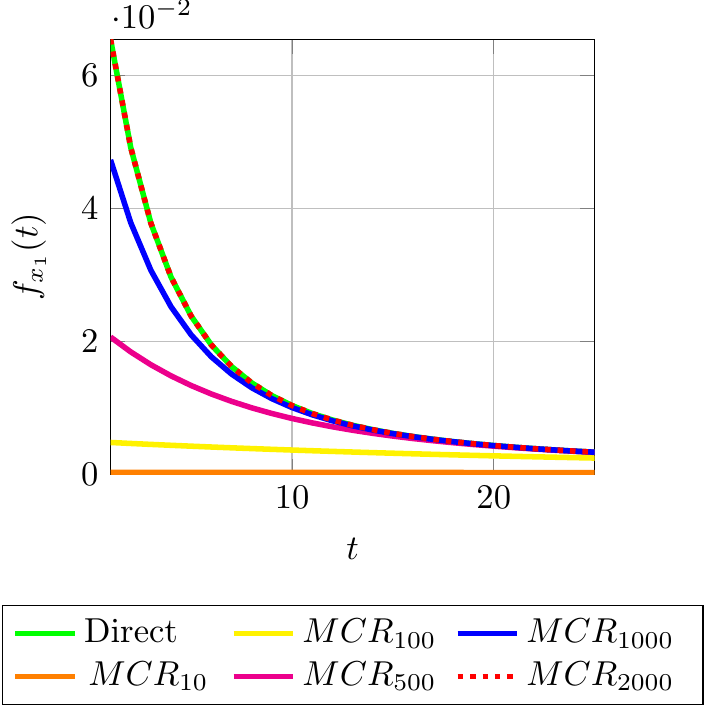}
\caption{Results on the dataset wolf by using the geometric heat equation for $t_M=25$. We compare the discrete feature
  descriptor $f_{x_1}(t)$ obtained by the direct solver and the MCR technique for different number of modes $N_{\max}$, shortened here as $MCR_{N_{\max}}$.
  Obviously, a small number of used modes (here $N_{\max}=10$ or $N_{\max}=100$) correlate to a slow rate of heat transfer. Therefore, the feature
  descriptors tend to almost constant functions when decreasing the number of modes, and this may lead to shape signature values that give
  not much information about the geometry of the given surface.}
\label{feature_descriptor_BE_MCR}
\end{figure}

After this brief discussion, let us now turn again to the MCR method.
The analogous representation as above is described in MCR which is based on dominant modes (low frequencies).
The solution of the reduced model $\dot{\mathbf{w}} (t) = \Lambda_r \mathbf{w} (t)$ of order $r$ is given by
\begin{align}
\mathbf{w}(t)=\underbrace{e^{\Lambda_r t}}_{\approx I}\mathbf{w}(0)\approx\mathbf{w}(0)
\label{solutionreducedHeatVector}
\end{align}
so that with $\mathbf{u}(t)\approx V_r \mathbf{w}(t)$ follows
\begin{align}
\mathbf{u}(t)\approx V_r \mathbf{w}(t)= V_r \mathbf{w}(0)= V_r V_r^\top D \mathbf{u}(0)
\label{solutionreducedHeatVector2}
\end{align}
Consequently, the pointwise MCR feature descriptor results in
\begin{align}
  f_{x_j}(t) & \approx \bigl( V_r V_r^\top D \mathbf{u}(0) \bigr)^\top
  \mathbf{\hat{e}_j}\nonumber\\
  &= \bigl( V_r V_r^\top \mathbf{\hat{e}_j} \bigr)^\top
    \mathbf{\hat{e}_j} = \left( V_r \mathbf{ y_{r,j}} \right)^\top \mathbf{\hat{e}_j} = c_{r,j}
\label{solutionreducedHeatVector3}
\end{align}
Having in mind the approximately constant behaviour of the approximative solution \eqref{solutionHeatVector4}
and comparing this with \eqref{solutionreducedHeatVector3} shows that the MCR descriptor exhibits
in general a slow rate of heat transfer.
Let us note that this shows that the MCR shape signature is not highly precise in terms of the indicated geometric location on a shape,
especially when using a small number of modes as shown for $f_{x_1}(t)$ in Figure \ref{feature_descriptor_BE_MCR}.


This observation requires an appropriate adaptation of the MCR method concerning the computation of the feature descriptor
in \eqref{solutionreducedHeatVector}-\eqref{solutionreducedHeatVector3}. An useful way to enhance without too much computational
effort the geometric information for the shape matching process is to modify the numerical signatures via an
\emph{adapted temporal domain}. As we will find, an \emph{adapted time} $t^\star$ with $t^\star \gg t_M$ for
small eigenvalues $0<\abs{\lambda_i}\ll 1$ causes $e^{\lambda_i t^\star} \ll e^{\lambda_i t_M}\approx 1$.
Let us explain, that we have $\lambda_i \leq 0$, so that the latter inequality bears the meaning that we obtain by rescaling time
in the indicated way a more significant spreading of values, so that the eigenvalues are more discriminative.
This improves the physical characteristics and consequently the diversity of the MCR signatures.
\par
Let us begin with the following thoughts. Obviously, the rate of heat transfer of the
reduced feature descriptor of order $r$ depends on the fastest mode ${\lambda_r}$. Moreover, a higher order reduced
model $r'>r$ corresponds to faster heat transfer, which in turn implies $t^\star_{r'}<t^\star_{r}$ to be suitable
in order to get a more discriminative descriptor value.
Based on these considerations, we will modify $[0,\mathit{t_M}]$ to $[0,\mathit{t^\star}]$ in the following way.

We propose to adapt the temporal domain in a simple way by the function
\begin{align}
t^\star(\lambda_r)=\frac{k}{\sqrt{\abs{\lambda_r}}}
\label{newTime}
\end{align}
where $k$ is a constant. This still to be determined parameter is defined in our construction
via the condition
\begin{align}
t_M\stackrel{!}{=}t^\star(\lambda_N)=\frac{k}{\sqrt{\abs{\lambda_N}}}
\label{newTimek}
\end{align}
which must be fulfilled for the reduced model of full order $N$.
Therefore, making use of \eqref{newTimek} the modified integration domain $[0,\mathit{t^\star(\lambda_r)}]$ is specified by 
\begin{align}
t^\star(\lambda_r)=\frac{t_M \sqrt{\abs{\lambda_N}}}{\sqrt{\abs{\lambda_r}}}
\label{newTimeHeat}
\end{align}
The corresponding time increment is easily calculated via $\tau=\frac{t^\star}{M}$.  
One may notice, that the square root prevents an extremely large time domain (eigenvalues $\abs{\lambda_r}$
grow exponentially) which would lead anew to less descriptive shape signatures.
\par
\begin{figure}[t]
\centering
\includegraphics[width=0.38\textwidth]{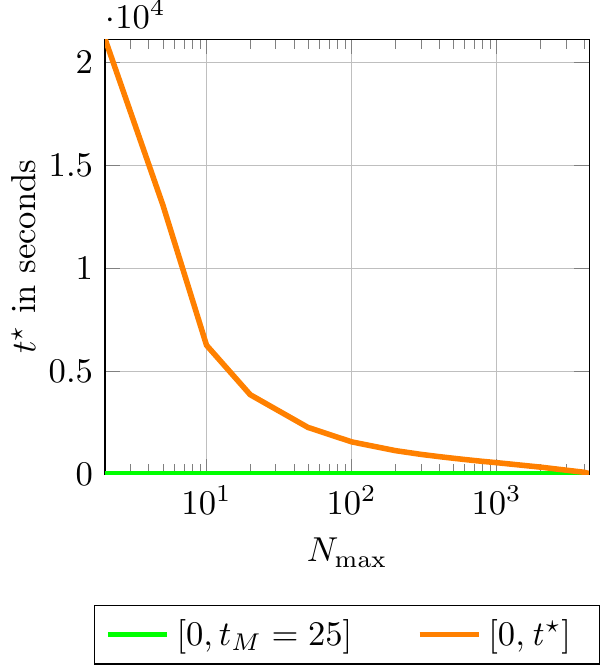}
\caption{Results on the dataset wolf by using the geometric heat equation. Visualisation of the
  \emph{original} $[0,\mathit{t_M}=25]$ and the \emph{adapted} $[0,\mathit{t^\star}]$ temporal
  domain using the formula \eqref{newTimeHeat} for a different number of modes $N_{\max}\in[2,4344]$.}
\label{modified_Time}
\end{figure}
\begin{figure}[t]
\centering
\includegraphics[width=0.48\textwidth]{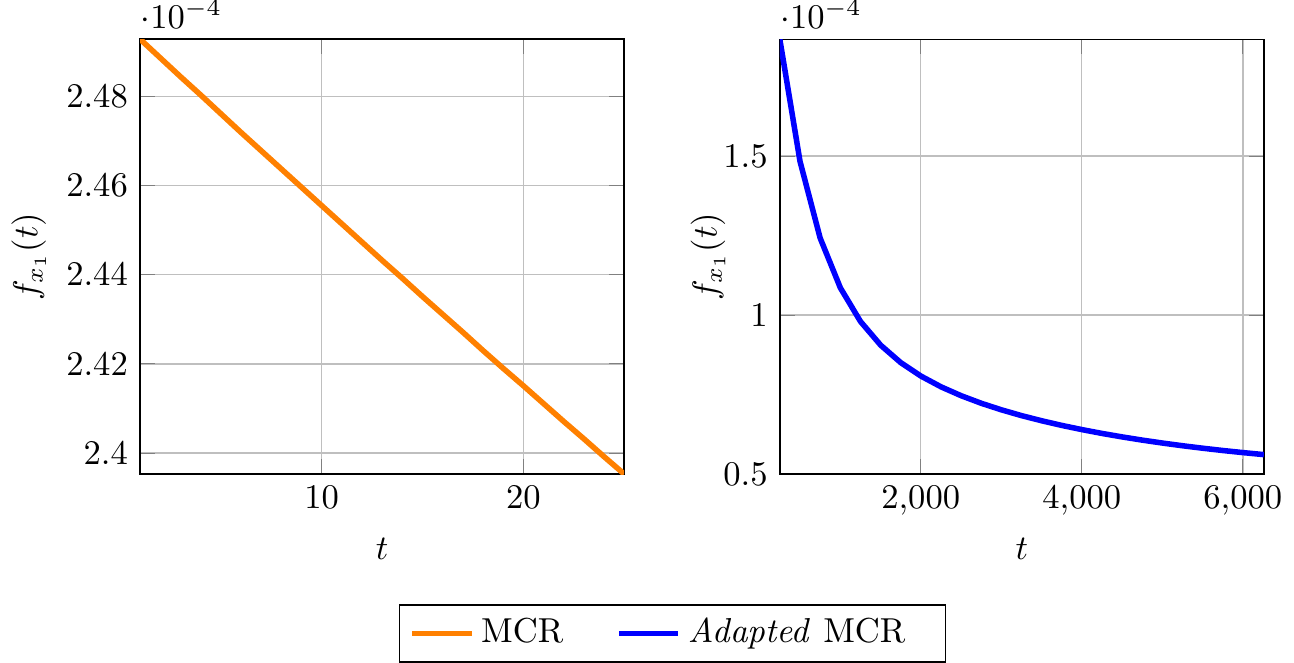}
\caption{Results on the dataset wolf by using the geometric heat equation for $M=25$ time levels. We compare the discrete
  feature descriptor $f_{x_1}(t)$ between the \emph{original} $[0,\mathit{t_M}=25]$ and the \emph{adapted} $[0,\mathit{t^\star}]$
  temporal domain using the formula \eqref{newTimeHeat} computed by the MCR technique for  $N_{\max}=10$. Note that the left
  figure shows a nearly constant function whereas a decaying exponential function as on the right after rescaling is the desired
  result. Clearly, the modified integration domain yields
  a more suitable numerical signature.} 
\label{feature_descriptor_MCR_10}
\end{figure}

As an example for illustrating the usefulness of our rescaling, the mentioned strategy yields for the wolf dataset
the modified integration domain $[0,\mathit{t^\star}]$ shown in Figure \ref{modified_Time}.
Based on this modification the adapted feature descriptor $f_{x_1}(t)$, illustrated in
Figure \ref{feature_descriptor_MCR_10}, demonstrates a more suitable numerical signature. \par
\begin{figure}[t]
\centering
\includegraphics[width=0.48\textwidth]{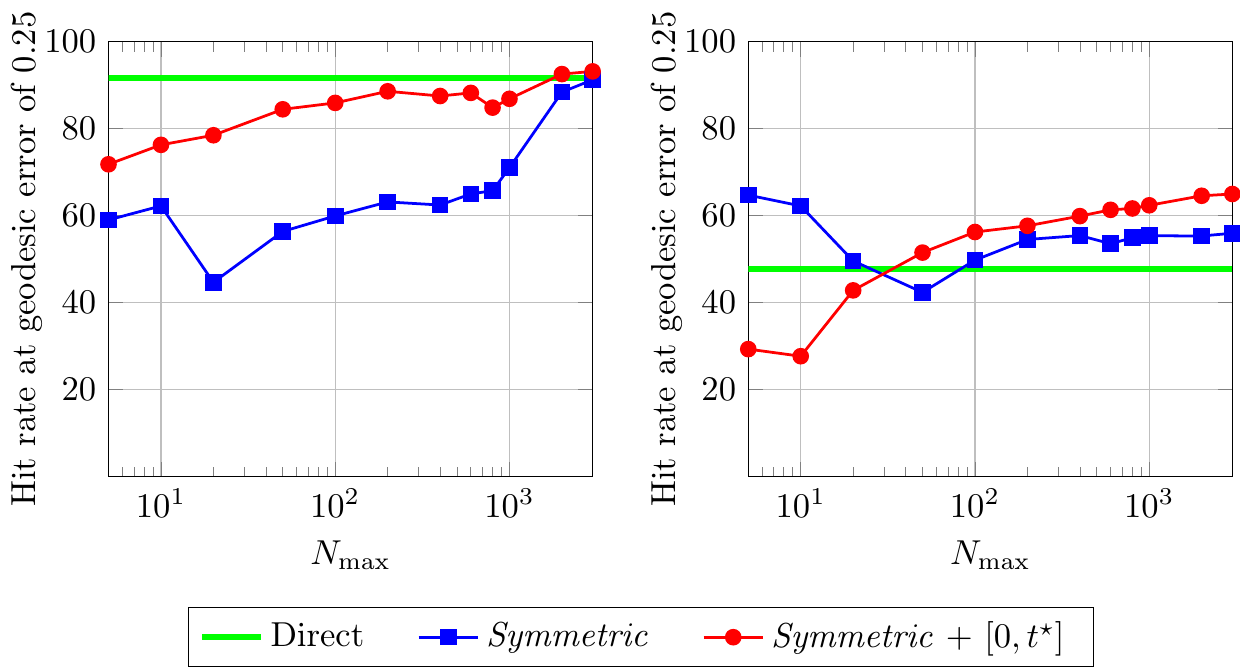}
\caption{Results on the dataset wolf (left) and baby (right) by using the geometric heat equation. We compare the geodesic error at $0.25$ between the \emph{original} $[0,\mathit{t_M}=25]$ and the \emph{adapted} $[0,\mathit{t^\star}]$ temporal domain using formula \eqref{newTimeHeat} computed by the MCR technique for different number of modes $N_{\max}\in[5,3000]$. The eigenvalues and eigenvectors are constructed by using the \emph{similarity} transformation.}
\label{geoerror_similar_T25_modified_time_MCR}
\end{figure}

Finally, we compare the performances of the MCR technique including the similarity transformation by using
the \emph{original} $[0,\mathit{t_M}=25]$ and the \emph{adapted} $[0,\mathit{t^\star}]$ temporal domain for the given
datasets. The corresponding visualisation of the evaluation is presented in
Figure \ref{geoerror_similar_T25_modified_time_MCR}.
In general, the substantial increase of the stopping time $t^\star(\lambda_r)$ of the reduced model improves the
geodesic error accuracy. However, the results may differ for a small spectrum
as in the case of the baby dataset. A possible heuristic explanation could be based on the nature of the first few
modes that correspond to a low pass filtered shape signature, meaning that one may obtain a good matching of a coarse
shape representation in some cases.
\par
In what follows, we will denote the solver including both abovementioned improvements as \emph{optimised} MCR method. 


\subsection{Implementation of the Optimised MCR Heat Signature}

We now describe briefly the simple implementation of the MCR signature which makes it a candidate for practical
shape analysis purposes.
\par
As already mentioned, the feature descriptors have to be computed for each point and on each shape. Nevertheless,
the computation of the MCR signatures can be done simultaneously for all locations $x_i\in\mathcal{M}$ on a given shape
by putting it into the format of a matrix-vector-multiplication, which additionally improves the computational efficiency.
\par
To this end, solving the geometric heat equation by \eqref{BEreducedNew} can be reformulated as
\begin{align}
W^{k} =&  PW^{k-1}, \quad W^{0}=V_r^\top DU^{0}= V_r^\top \label{heatImplementation}
\end{align}
with $V_r \in \mathbb{R}^{N\times r}$,
$W\in\mathbb{R}^{r\times N}$ and $U^0=D^{-1}$ due to $\mathbf{u_i}(0)=(0,\dots,0,|\Omega_i|^{-1},0,\dots 0)^\top$.
The algorithm for computing the optimised MCR heat signature is described in Figure \ref{algorithmHeat}.
Let us note that the extraction of the feature descriptors $f_{x_i}$ in step 3d) of the algorithm is based
on the Hadamard product. Finally, the correspondence quality can be evaluated based on the selected measure.
\par
Let us mention, incidentally, the computation of the temporal domain $[0,\mathit{t^\star}]$ and the time increment $\tau$ in step 2 of the described algorithm is only done for the reference shape. This ensures to perform a correct matching by comparing the feature descriptors for different locations $x_i\in \mathcal{M}$ and $\widetilde{x}_j\in \widetilde{\mathcal{M}}$ on the same scale.
\begin{figure}
\begin{algorithm}[H]
\caption{Computation of MCR Heat Signature}
\begin{algorithmic}
\REQUIRE Matrices $D^{-1},W$, reduction parameter $r$, amount of time levels $M$, stopping time $t_M$ 
\ENSURE $\mathbf{f_{x_i}}$\\\hrulefill\\\vspace{0.1cm}
\textbf{Reference Shape}\\\hrulefill\vspace{0.1cm}
\STATE 1. Computation of $\Lambda_r \in\mathbb{R}^{r\times r}$  and $V_r \in \mathbb{R}^{N\times r}$\vspace{0.3cm}
\STATE a) Compute $r$ dominant eigenpairs $(\lambda,\mathbf{v})$ of \\
$B=D^{-\frac{1}{2}} WD^{-\frac{1}{2}}$ by
solving $B\mathbf{v}=\lambda \mathbf{v}$\vspace{0.1cm}
\STATE b) Calculation of $(\lambda, \mathbf{\widetilde{v}})$ of $L$ with $\mathbf{\widetilde{v}}=D^{-\frac{1}{2}}\mathbf{v}$
\\\hrulefill\vspace{0.1cm}
\STATE 2.) Computation of $[0,\mathit{t^\star}]$ using \eqref{newTimeHeat} \vspace{0.3cm}
\STATE a) Compute fastest mode $\lambda_N$ by solving $B\mathbf{v}=\lambda \mathbf{v}$ \vspace{0.1cm}
\STATE b) Calculate $t^\star(\lambda_r)=\displaystyle{\frac{t_M \sqrt{\abs{\lambda_N}}}{\sqrt{\abs{\lambda_r}}}}$, $\tau=\displaystyle{\frac{t^\star}{M}}$
\\\hrulefill\vspace{0.1cm}
\STATE 3.) Solving reduced system \eqref{heatImplementation}\vspace{0.3cm}
\STATE a) Set $W^{0}=V_r^\top$\vspace{0.1cm}
\STATE b) Compute $P=(I- \tau \Lambda_r)^{-1}$\vspace{0.1cm}
\STATE c) Solve $W^k=PW^{k-1}$ for $k=1,\ldots,M$\vspace{0.1cm}
\STATE d) Extract $\mathbf{f_{x_i}}$ by $f_{x_i}(t_k)=\sum\limits_{j=1}^r[(V_r\circ (W^k)^\top)]_{ij}$ \\
for $k=1,\ldots,M$\\\hrulefill\\\vspace{0.1cm}
\textbf{Transformed Shape}\\\hrulefill\vspace{0.1cm}\\ 
Repeat computations of Steps 1. and 3. 
\end{algorithmic}
\end{algorithm}
\caption{Algorithm for computing the \emph{optimised} MCR signature by solving the geometric heat equation.
  The procedure is analogously applicable to the geometric wave equation.}
\label{algorithmHeat}
\end{figure}
\subsection{Optimised MCR Wave Signature}
\label{MCRwaveSignature}
\begin{figure}[t]
\centering
\includegraphics[width=0.48\textwidth]{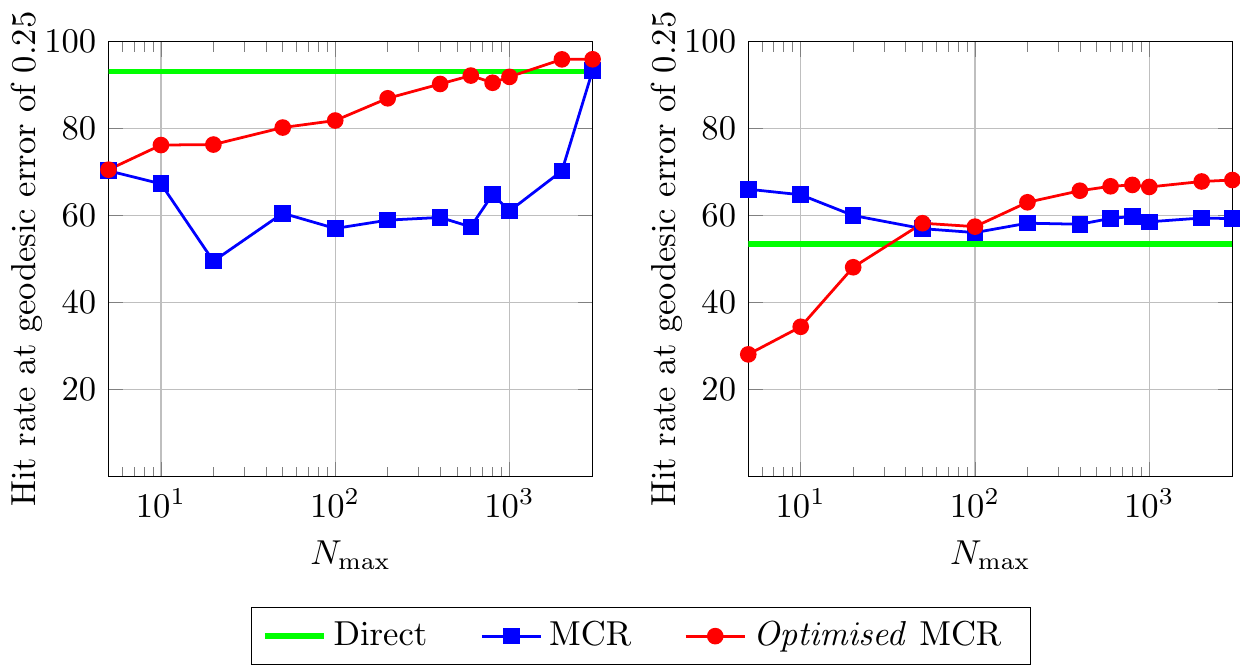}
\caption{Results on the dataset wolf (left) and baby (right) by using the geometric wave equation. We compare the geodesic error at $0.25$ between the \emph{original} $[0,\mathit{t_M}=25]$ and the \emph{adapted} $[0,\mathit{t^\star}]$ temporal domain using formula \eqref{newTimeWave} computed by the MCR technique for different number of modes $N_{\max}\in[5,3000]$. The eigenvalues and eigenvectors of the \emph{original} and \emph{optimised} MCR are constructed by solving the \emph{generalised} and \emph{symmetric} eigenvalue problem, respectively.}
\label{geoerror_MCR_wave}
\end{figure}

\begin{figure}[t]
\centering
\includegraphics[width=0.48\textwidth]{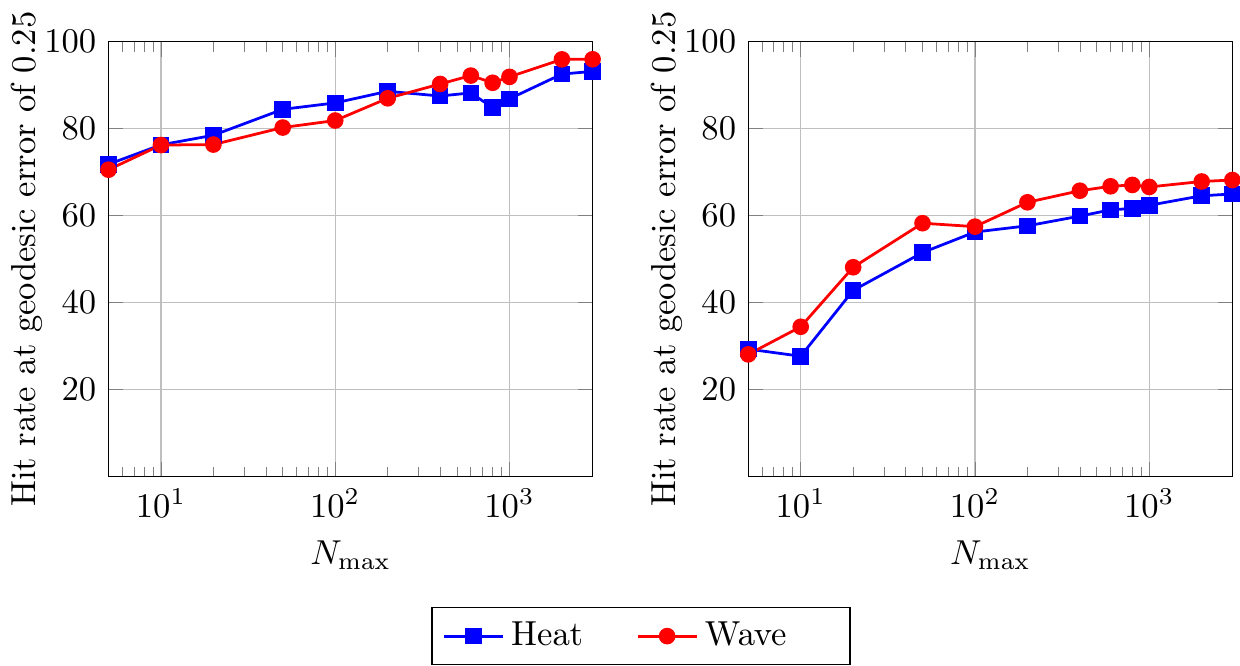}
\caption{Results on the dataset wolf (left) and baby (right) by applying the \emph{optimised} MCR technique. We compare the geodesic error at $0.25$ between using the geometric heat and wave equation for different number of modes $N_{\max}\in[5,3000]$.}
\label{geoerror_comparison_MCR_heat_wave}
\end{figure}

The abovementioned aspects can also be used to solve the geometric wave equation numerically by the MCR technique.
However, based on our experiments we propose to define the modified integration domain $[0,\mathit{t^\star(\lambda_r)}]$
in the following way 
\begin{align}
t^\star(\lambda_r)=\frac{t_M \sqrt[4]{\abs{\lambda_N}}}{\sqrt[4]{\abs{\lambda_r}}}
\label{newTimeWave}
\end{align}
The implementation (cf. Figure \ref{algorithmHeat}) for computing the optimised MCR wave signature remains identical
to the algorithm presented above, except for the calculation of $t^\star(\lambda_r)$ and for solving the reduced system in Step 2b) and Step 3, respectively. The modifications in Step 3 are the following: 
\begin{itemize}
\item[b)] Compute $P=(I- \tau^2 \Lambda_r)^{-1}$\\
\item[c)] Solve $W^1=PW^0$ and $W^k=P(2W^{k-1}-W^{k-2})$ for $k=2,\ldots,M$\vspace{0.1cm}
\end{itemize}
For completeness the results of the standard and optimised MCR wave signatures for the given datasets are presented in
Figure \ref{geoerror_MCR_wave}. Moreover, a comparison of the optimised MCR method by using the geometric heat and wave equation
is illustrated in Figure \ref{geoerror_comparison_MCR_heat_wave}.

\section{Reference Model: Spectral Methods}
The spectral methods of the heat kernel signature (HKS) and the wave kernel signature (WKS) that we employ for comparison
with our approach are based on the geometric heat equation \cite{SOG} and the geometric Schr\"odinger equation \cite{ASC}, respectively.
Let us note here, that also the process described in \cite{MORCC} is technically related to our appraoch as it is a
time evolution method and based on diffusion. However, as already noted in \cite{MORCC}, the technique is computationally
very demanding compared to HKS and WKS construction. Thus we refrain from comparing to this method since one of our main
aims is computational efficiency.

For constructing spectral descriptors, it is assumed that the solution of both equations will take the 
form $u(x,t)= \phi(x) \theta(t)$ due to the fact that the underlying PDEs are linear and homogeneous. 
This approach works because if the product of two functions 
$\phi$ and $\theta$ of independent variables $x$ and $t$ is a constant, each function must separately be a constant. 
Therefore, one may separate the equations to get a function of only $t$ and $x$, respectively: 
\begin{align}
\kappa \dfrac{ \partial_t \theta(t)}{\theta(t)} =  \dfrac{\Delta_\mathcal{M} \phi(x)}{\phi(x)} =\text{const} =-\lambda 
\end{align}
where $\kappa$ summarises both equations ($\kappa=1$ for geometric heat equation, $\kappa=i$ for geometric Schr\"odinger equation), 
and  $-\lambda$ is called the separation constant which is arbitrary for the moment.
This leaves us with two new equations, namely an ODE for the temporal component 
\begin{align}
\partial_t \theta(t)  =   - \kappa \lambda \theta(t), \quad t\in [0,T]
\end{align}
and the spatial part takes the form of the Helmholtz equation
\begin{align}
{\begin{cases} \Delta_\mathcal{M} \phi(x)  = - \lambda \phi(x),
&{\mbox{}} x\in \mathcal{M}
\\   
\langle \nabla_\mathcal{M} \phi , n \rangle = 0,
&{\mbox{}} x\in \partial \mathcal{M}
\end{cases}} 
\label{helmholtz-1}
\end{align}
where the constant $\lambda$ has here the meaning of the operator's eigenvalue.
\par
The Laplace-Beltrami operator is a self-adjoint operator on the space $L_2 (\mathcal{M})$
(since we assumed the shapes to be compact). 
This implies that the Helmholtz equation
has an infinite number of non-trivial solutions $\phi_k$ for certain eigenvalues $\lambda_k$
and corresponding eigenfunctions, which is a result of the spectral theorem \cite{Sauvigny2012}. 
Consequently, \eqref{helmholtz-1} takes on the format
\begin{align}
{\begin{cases} \Delta_\mathcal{M} \phi_k(x)  = - \lambda_k \phi_k(x),
&{\mbox{}} x\in \mathcal{M}
\\   
\langle \nabla_\mathcal{M} \phi_k , n \rangle = 0,
&{\mbox{}} x\in \partial \mathcal{M}
\end{cases}} \quad k=1,2,\ldots
\end{align}
Concerning the ordered spectrum of eigenvalues
$0=\lambda_1 < \lambda_2 \leq \ldots$ and the corresponding eigenfunctions $\phi_1 , \phi_2 , \ldots $
let us note that the latter form an orthonormal basis of $L_2 (\mathcal{M})$. Moreover, in case of Neumann
boundary conditions (or no shape boundaries) constant functions are solutions of the Helmholtz equation for the
first eigenvalue with zero value.
\par
It is well known that the eigenfunctions are a natural generalisation of the classical Fourier basis for working with 
functions on shapes. Let us also note that the physical interpretation of the Helmholtz equation is the following. 
The shape of a 3D object can be thought of as a vibrating membrane, the $\phi_k$ can be interpreted as its vibration modes 
whereas $\lambda_k$ have the meaning of the corresponding vibration frequencies, sorted from low to high. \par
For each index $k$ one thus obtains an ODE \eqref{ode} which can be solved
by integration using the indefinite integral:
\begin{align}
\int \dfrac{\mathrm{d}\theta(t)}{\theta(t)}    =   - \int \kappa \lambda_k \mathrm{d}t 
\quad \Longrightarrow \quad
\theta(t)= \alpha_k e^{- \kappa \lambda_k t }
\end{align}
where the integration constant $\alpha_k$ should satisfy the initial condition of the $k$-th eigenfunction.
The final product solution then reads as
\begin{align}
u_k(x,t)= \alpha_k e^{- \kappa \lambda_k t } \phi_k(x) \, 
\end{align}
The principle of superposition says that if we have several solutions to a linear homogeneous differential equation 
then their sum is also a solution. Therefore, a closed-form solution of the geometric
heat equation in terms of a series expression can be written as
\begin{align}
u(x,t)= \sum_{k=1}^{\infty} \alpha_k e^{- \lambda_k t } \phi_k(x) \, 
\end{align}
and the solution of the geometric Schrödinger equation reads as
\begin{align}
u(x,t)= \sum_{k=1}^{\infty} \alpha_k e^{- i\lambda_k t } \phi_k(x) \, 
\end{align}
where the coefficients $\alpha_k$ are chosen to fulfil the initial condition. 

\paragraph{Heat Kernel Signature.} The coefficients $\alpha_k$ in our expansion 
can be computed by using  the $L_2$ inner product
\begin{align}
\alpha_k= \langle u_0, \phi_k \rangle_{L_2(\mathcal{M})} = \int\limits_{\mathcal{M}} u_0(y) \phi_k(y)  \, \mathrm{d}y 
\end{align}
such that one may compute
\begin{align}
u(x,t)&= \sum_{k=1}^{\infty} \Bigg( \int\limits_{\mathcal{M}} u_0(y) \phi_k(y)  \, \mathrm{d}y \Bigg) e^{- \lambda_k t } \phi_k(x)\\
&=  \int\limits_{\mathcal{M}} u_0(y) \underbrace{\bigg( \sum_{k=1}^{\infty}   e^{- \lambda_k t } \phi_k(y)  \phi_k(x) \bigg)}_{=K(x,y,t)} \, \mathrm{d}y \, 
\end{align}
The {\em heat kernel} $K(x,y,t)$ describes the amount of heat transmitted from $x$ to $y$ after time
$t$. By setting the initial condition to be a delta heat distribution with $u_0(y)=\delta_x(y)$ at the position $y$, we thus obtain
according to \cite{SOG} the heat kernel signature
\begin{align}
\text{HKS}(x,t) = 
\sum_{k=1}^{\infty}    e^{- \lambda_k t }  \left| \phi_k(x) \right|^2 \, 
\end{align}
where the shifting property of the delta distribution $f(x)=\int_\mathcal{M} f(y) \delta_x(y) \, \mathrm{d}y$ was used. In accordance, the quantity
$\text{HKS}(x,t)$ describes the amount of heat present at point $x$ at time $t$.

\paragraph{Wave Kernel Signature.} The WKS \cite{ASC} is defined to be the time-averaged probability of detecting a particle of a certain energy distribution at the point $x$, formulated as 
\begin{align}
\text{WKS}(x,t) &= \lim\limits_{T\rightarrow \infty} \dfrac{1}{T} \int\limits_{0}^{T} |u|^2 \mathrm{d}t 
 = \sum\limits_{k=1}^{\infty}  \left| \alpha_k  \right|^2 \left| \phi_k(x) \right|^2 \, 
\end{align}
Furthermore, $\alpha_k=\alpha(e_k)$ becomes a function of the energy distribution $e_k$ of the quantum mechanical
particle and can be chosen as a log-normal distribution, i.e.
\begin{align}
 \left| \alpha_k  \right|^2=   \exp \bigg( \frac{-(e-\log \lambda_k)^2}{2 \sigma^2} \bigg) \, 
\end{align}
where the variance of the energy distribution is denoted by $\sigma$, see again \cite{ASC} for more details.

\par

The eigenfunctions and eigenvalues of the discrete Laplace-Beltrami operator are computed by performing the generalised 
eigendecomposition
\begin{align}
L \phi_k  = -\lambda_k  \phi_k \, ,  \quad k=1,\ldots,N \, 
\end{align}

\section{Experimental Results: Optimised MCR vs. Kernel-based Methods}
In the following, we give a detailed quantitative evaluation of the proposed MCR technique compared to the kernel-based methods.
In the total, we will compare the developed optimised MCR techniques relying on both geometric heat and wave equation, respectively,
with HKS and WKS. To this end, we benchmark the methods at hand of the complete TOSCA dataset whose shapes are almost isometric.
\par
Let us note again, that the comparison of the optimised MCR wave signature and WKS is not based on the same PDE model although the
latter notion (wave kernel signature) suggests this, as the latter method is based on the Schrödinger equation. Concerning this point,
let us note
first that a recent work \cite{DBH-2017} has shown that the numerical descriptor based on the geometric wave equation gains better
results than the numerical descriptor based on the geometric Schrödinger equation. Secondly, within the class of kernel-based methods
the WKS can be considered as a competitive descriptor. Therefore, we compare the optimised MCR wave signature in its relation to the
WKS due to their distinct role in their respective class.

\paragraph{Evaluation Measure.}
The correspondence quality will again be measured by using the  \emph{geodesic error}.

\paragraph{Technical Remarks.}
The TOSCA dataset we investigate includes several classes of almost isometric shapes. In detail, it contains 76 shapes
divided into 8 classes (humans and animals) of varying resolution (3K to 50K vertices).
Furthermore, for some introductory experiments the already used wolf and the baby dataset are evaluated.
\par
The experimental comparison basically considers the implementations of the kernel-based methods and
its parameter settings as described in \cite{ASC,SOG}.
For all methods we compute the feature descriptors sampled at 100 points. On this basis, the adapted temporal
domain $[0,\mathit{t^\star}]$ and the uniform time increment $\tau$ are calculated by using $t_M=25$ and $M=100$.
\par
All experiments were done in MATLAB R2018b with an Intel Xeon(R) CPU E5-2609 v3 CPU. The eigenvalues and eigenvectors
are computed by the Matlab internal function \textit{eigs}.

\subsection{Evaluation of the Geodesic Error}\label{EvaluationonGeodesicError}
In the following, the evaluation is subdivided into two parts. First, we compare the correspondence quality of the
optimised MCR technique and the kernel-based methods based on the given datasets wolf and baby.
Subsequently, all methods are benchmarked at hand of the complete TOSCA dataset. Let us note in this context,
that the numerical advances proposed here enable for the first time the evaluation of the complete TOSCA dataset,
as documented in this paper. As indicated beforehand, this together with the proposed technical improvements
represents also a significant experimental extension with respect to \cite{BDB-2018}.
\paragraph{Evaluation on Wolf and Baby Datasets}
\begin{figure}[t]
\centering
\includegraphics[width=0.48\textwidth]{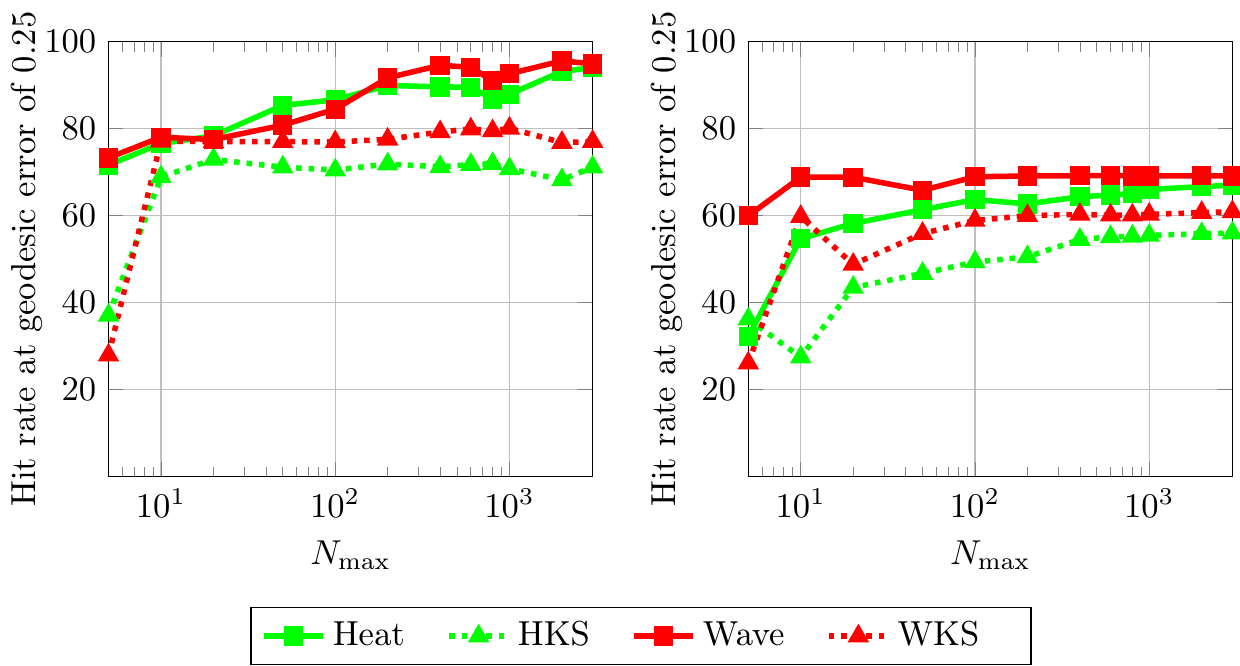}\\
\caption{Results on the dataset wolf (left) and baby (right). Comparison of the geodesic error at $0.25$ between
  the \emph{optimised} MCR technique and the kernel-based methods for different number of modes $N_{\max}\in[5,3000]$.
  The eigenvalues and eigenvectors for MCR and the kernel-based methods are computed by solving the \emph{symmetric}
  and the \emph{generalised} eigenvalue problem, respectively.}
\label{Wolf_Baby_Comparison}
\end{figure}
The results for both datasets, shown in Figure \ref{Wolf_Baby_Comparison}, clearly demonstrate a higher matching performance
by using the optimised MCR technique. In almost all cases, the MCR signatures outperform their competitive methods HKS and WKS.
\par
\begin{figure}[t]
\centering
\includegraphics[width=0.38\textwidth]{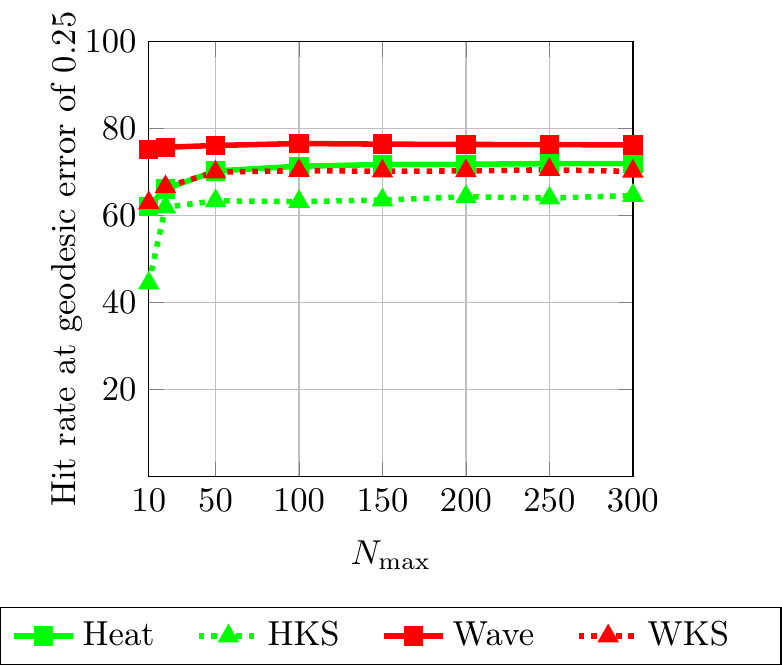}\\
\caption{Results on the TOSCA dataset. Comparison of the geodesic error at $0.25$ between the optimised MCR technique and
  the kernel-based methods for a varying number of modes $N_{\max}$. The eigenvalues and eigenvectors for MCR and the
  kernel-based methods are computed by solving the \emph{symmetric} and the \emph{generalised} eigenvalue problem, respectively.}
\label{Tosca_Comparison_geoerror}
\end{figure}

An interesting aspect is that the MCR heat signature yields slightly better results as WKS. Moreover, the depicted curves in
Figure \ref{Wolf_Baby_Comparison} show a saturation behaviour with respect to the number of used modes. Specifically, that
means that the correspondence quality can no longer be significantly improved after a certain number of used modes.
\paragraph{Evaluation on TOSCA Dataset}
First, let us also discuss here the number of the ordered modes for the mentioned methods. In many publications on the
kernel-based methods this number is set manually to a fixed value (e.g.\ $N_{\max}=300$ is often used). Because of its impact
on computational efficacy, we were motivated to consider here the examined measurements and the amount of used eigenvalues
for the TOSCA dataset, see the results illustrated in Figure \ref{Tosca_Comparison_geoerror}.
\par
Over the whole dataset, we obtain qualitatively identical results as for the examples wolf and baby. The proposed MCR technique
outperforms the kernel-based methods in terms of the geodesic error. For both PDEs, the optimised MCR method provides 
around 5-10\% higher correspondence quality than the corresponding kernel-based approach. Again, in general the MCR heat signature
slightly exceeds the performance of WKS. Additionally, the experiment indicates that for all methods a saturation behaviour is achieved at a
small spectrum of $N\approx 100$ modes.

\subsection{Soft Correspondence}
\label{SoftCorrespondence}

For two discrete shapes $\mathcal{M}_d$ and $\mathcal{\widetilde{M}}_d$ with underlying point clouds 
$P=\{x_1,\ldots,x_N \}$ and $\widetilde{P}=\{\widetilde{x}_1,\ldots, \widetilde{x}_{\widetilde N} \}$,
the matrix representation of the correspondence map $S_{\{0,1\}}: \widetilde{P}\times P \rightarrow \{0,1\}$ can be encoded by 
\begin{equation}
{S_{\{0,1\}}}_{ji}=
 \begin{cases}
 1, & {\mbox{if }} d_f(x_i,\widetilde{x}_j) \leq d_f(x_i , \widetilde{x}_k), \quad k= 1, \ldots,\widetilde N \\
 0, & \text{else}
 \end{cases}\label{PtoP}
\end{equation}
where $S_{\{0,1\}}$ is a $(\widetilde{N} \times N)$ binary assignment matrix. For finding a single corresponding counterpart
e.g.\ of $x_i \in P$ on $\mathcal{\widetilde{M}}_d$,
we construct a vector-valued indicator function $h: P \rightarrow \{0,1\}$ with
\begin{align}
 \mathbf{h_i}(x_k) = 
 {\begin{cases} 1         ,
&{\mbox{if }} i=k
\\  
0,
&{\mbox{else }} 
\end{cases}}
\end{align} 
By performing matrix-vector multiplication $\mathbf{\widetilde{h}_j}={S_{\{0,1\}}} \mathbf{h_i}$,
we obtain the indicator vector $\mathbf{\widetilde{h}_j}$ for the case the
tuple $(x_i, \widetilde{x}_j)\in P \times \widetilde{P}$ is a corresponding pair. \par
However, due to strong elastic deformations, noisy shapes or intrinsic symmetries (i.e.\ inherent ambiguities)
the construction \eqref{PtoP} may generate misleading matchings. Therefore, the correspondence map $S_{\{0,1\}}$ is in
practice neither injective nor surjective, even in the case $N = \widetilde{N}$.
\par
By allowing values between 0 and 1, the correspondence map can be cast as a
soft correspondence map $S_{[0,1]}: \widetilde{P}\times P  \rightarrow [0,1]$, as used for example in \cite{Eisenberger2018}, 
which can be interpreted as a correspondence probability by using normalisation 
\begin{align}
 \sum\limits_{i=1}^{N} {S_{[0,1]}}_{ji} =1, \quad j=1,\ldots , \widetilde{N}
\end{align}
For a certain point $x_i\in P$ on the reference shape the map encodes the transition probabilities
to the points on $\widetilde{P}$ based on the feature distances, whereby high probabilities correspond
to low feature distances and vice versa. 
The arising soft correspondence map allows us to express the probability transition
of all points from a reference shape onto a target shape.
\par
\begin{figure*}
	\centering
	\begin{tikzpicture}[]
	\node[ ] at (0,0) {\includegraphics [width=1\linewidth]{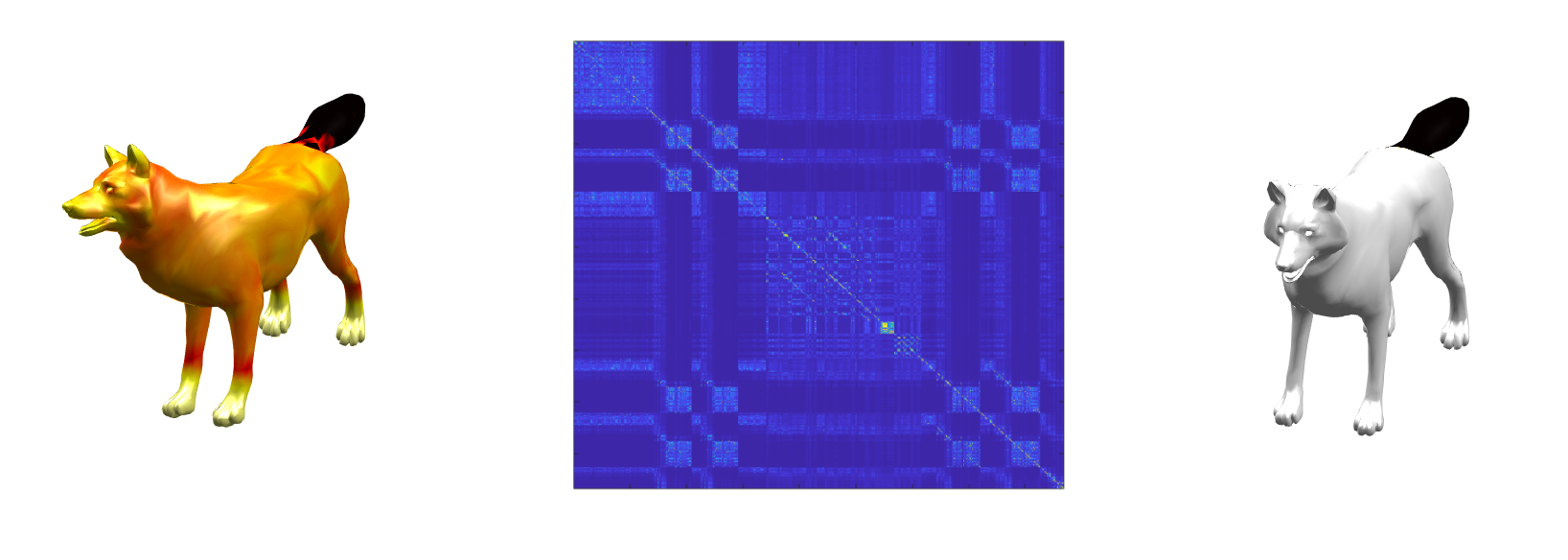}};
	\node[text width=3cm,scale=2.5] at (7.75,0) {$\times$};
	\node[text width=3cm,scale=2.5] at (-0.2,0) {$=$};
		\node[text width=3cm,scale=1.5] at (9.0,-1.7) {$\mathcal{M}$};
	\node[text width=3cm,scale=1.5] at (-5.4,-1.2) {$\mathcal{\widetilde{M}}$};
	\end{tikzpicture}
	\caption{Idea of mapping indicator functions: An indicator function (visualised by a colour map on the shape)
          $ \mathbf{h}\in \mathcal{M}$ is defined such that it is one at the tail of the dog and zero else. Using the soft
          correspondence matrix, we define the probability of the matched indicator function $\mathbf{\widetilde{h}}\in \widetilde{\mathcal{M}}$
          by performing $\mathbf{\widetilde{h}}=S_{[0,1]} \mathbf{h}$. For demonstration the optimised MCR heat method with $300$ modes was used.
          The probability of matching the indicator function on the tail of the transformed wolf is very high, while the
          colour encodes the probability ranging from white (almost zero probability) to black (probability of $10\%$ and higher). Since the ground truth is the identical labeling $(i,i)$, the soft correspondence matrix has a diagonal dominant structure.}
	\label{wolf1}
\end{figure*}

Let us briefly sketch an example for the usefulness of the approach. For this we make use of Figure \ref{wolf1}
where we consider in particular the tail of the wolf. Similar to the pointwise case, we construct a discrete indicator
function $\mathbf{h}$ with
\begin{align}
 \mathbf{h}(x_i) = 
 {\begin{cases} 1,
&{\mbox{if }} x_i \in \text{Tail} \subset \mathcal{M}
\\  
0,
&{\mbox{else }} 
\end{cases}}
\end{align} 
By performing $\mathbf{\widetilde{h}}=S_{[0,1]} \mathbf{h}$, the entries of $\mathbf{\widetilde{h}}$ contain the probability of the
matched indicator function on the transformed shape.

\par
\begin{figure}
	\centering
	\begin{tikzpicture}[]
	\node[ ] at (0,0) {\includegraphics [width=0.95\linewidth]{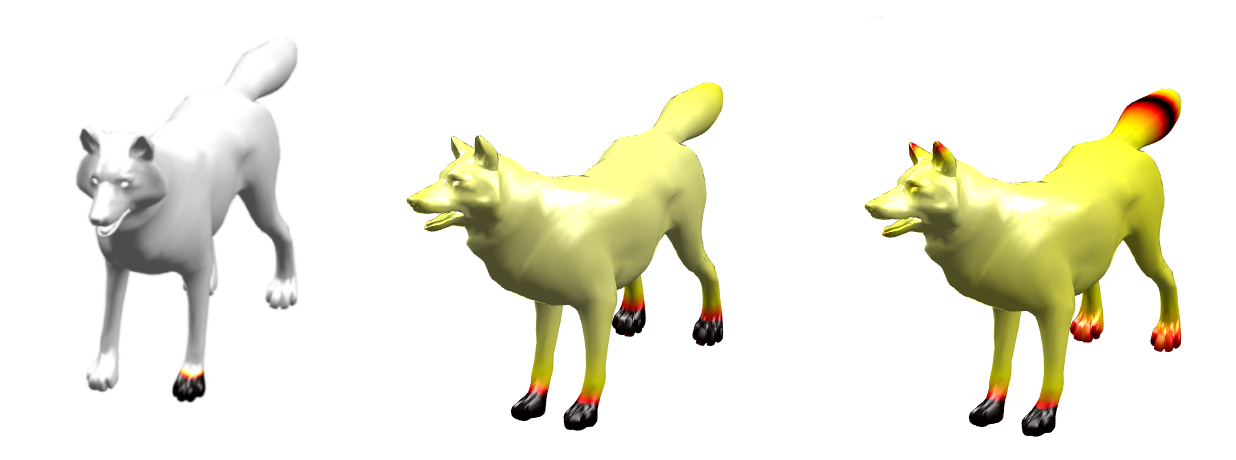}};
	\node[text width=3cm] at (0.4,1.4) {MOR Heat};
	\node[text width=3cm] at (3.8,1.4) {HKS};
	\end{tikzpicture}
	\caption{Mapping of the indicator function defined on a certain region of the wolf shape (the front paw)
          to a transformed version of it.
          The colour encodes the probability of the matching, being black for a high probability and white for
          almost zero probability. The optimised MCR heat approach leads in this example to a better interpretation
          compared to the HKS, detecting similar regions of the shape reliably. }
	\label{wolf2}
\end{figure}

Using soft correspondence maps in combination with the MCR technique may yield beneficial results
when addressing problems like correspondence of shape segments or detection of intrinsically
similar regions. The intuition behind this is that intrinsic feature descriptors have almost identical, low feature distances for
whole regions where the intrinsic geometry is similar, and that such region-based information could favourably be described by 
low-frequency modes as used within MCR, whereas small details are described by high-frequency modes; see
Figure \ref{wolf2} for an example.

\subsection{Sparsity of the Soft Correspondence Map}

The soft correspondence map $S_{[0,1]}$ is a dense $(\widetilde{N} \times N)$-matrix, being increasingly
cumbersome for large $\widetilde{N}$ and $N$. Especially for large shapes, the computation, storage and manipulation of
$S_{[0,1]}$ might be computationally expensive.
Therefore, it would be desirable to reduce the information contained in the dense matrix $S_{[0,1]}$,
ideally computing from it the binary assignment matrix $S_{\{0,1\}}$ , however this task is in practice not trivial.
\par
In the following, we demonstrate experimentally that it is possible to give $S_{[0,1]}$ a
sparse structure while keeping its useful meaning for pointwise shape correspondence.
The sparse structure appears to be a reasonable and computationally efficient compromise between
the dense soft correspondence matrix and the ideal binary assignment matrix.
\par
To increase sparsity we set iteratively small entries of $S_{[0,1]}$ (carrying a low probability) in ascending order
w.r.t.\ the size of entries to zero. We assume in this context that the feature descriptor based on MCR delivers correct
corresponding pairs for low feature distances (high probability entries of $S_{[0,1]}$) only.
\par
Let us note that in each sparsification step one could normalise again the resulting matrix $S_{[0,1]}$
(which we do not give a different notion here as during sparsification it keeps its role as a soft correspondence
matrix) so that one could again interprete it as probability transition matrix.

\paragraph{Evaluation Measure}
In order to measure sparsity of the soft correspondence matrix, we define the \textit{density} of $S_{[0,1]}$ by taking the
ratio of the numbers of non-zero elements to the total number of elements
(thus e.g.\ relative density of 100\% means we have a full matrix).
When using the ground truth in our experiments, we introduce the \textit{soft hit rate} by evaluating if
non-zero elements are transition probabilities leading to correct correspondences. It is obvious that the
soft hit rate is correlated to the density of $S_{[0,1]}$, since the more non-zero elements are set to zero
the more correct correspondences might accidentally be removed by the sparsification routine.
\par
During sparsification, we seek a threshold density related to the individual example, such that
the resulting matrices $S_{[0,1]}$ still contain the complete pointwise information of the underlying,
ideal binary assignment matrix $S_{\{0,1\}}$. Therefore, we define the
\textit{minimum density} by taking the sparsest computed version of $S_{[0,1]}$ such that the soft hit
rate still yields $100\%$ correct correspondences. 
\paragraph{Evaluation on Wolf and Baby Datasets}
\begin{figure}[t!]
\centering
\includegraphics[width=0.48\textwidth]{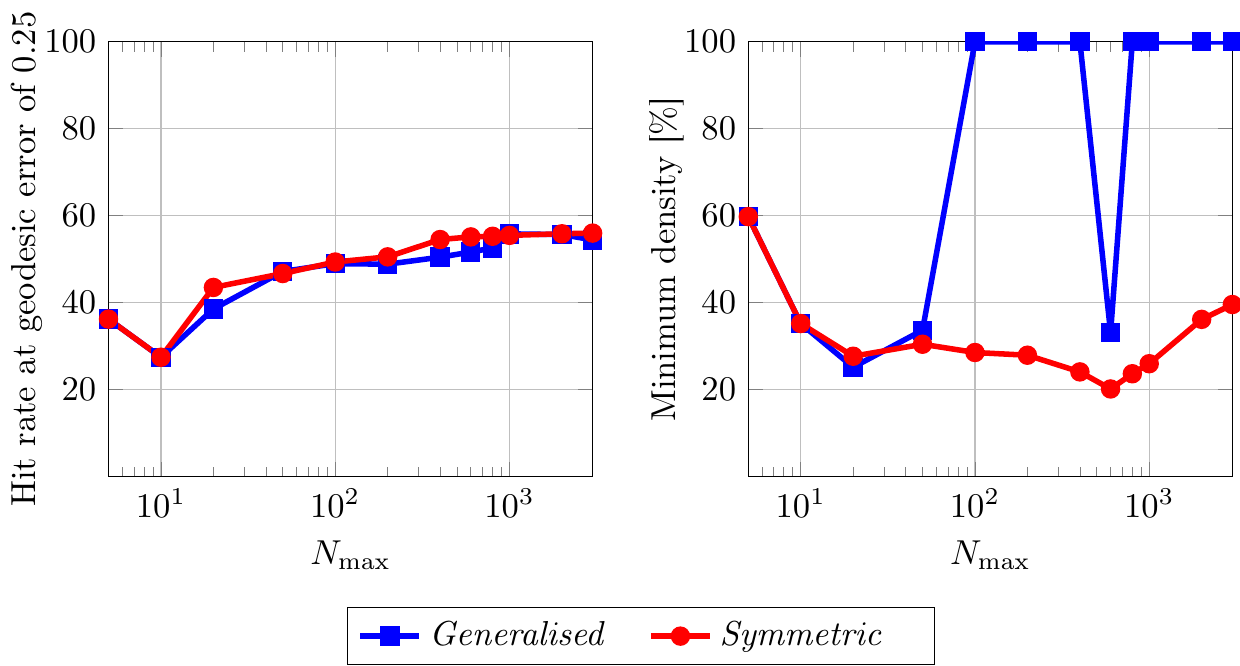}
\caption{Results on the dataset baby by applying the \emph{heat kernel signature} for a varying number of modes
  $N_{\max}\in[5,3000]$. Comparison of the geodesic error at $0.25$ (left) and the minimum density (right)
  between the construction of eigenvalues and eigenvectors by solving the generalised and the symmetric
  eigenvalue problem. Using the similarity transformation dramatically improves the sparsity of $S_{[0,1]}$. }
\label{HKS_generalized_similar}
\end{figure}
Before we go into a more detailed evaluation on the sparsity of $S_{[0,1]}$, we revise the issue of the eigenvalue
computation for the kernel-based methods. As demonstrated by the example of the baby dataset, cf.\ Figure \ref{HKS_generalized_similar},
the way of construction of eigenvalues and eigenvectors affects the matching results.
Even if the geodesic error accuracy remains almost identical, the modes computed by solving the generalised eigenvalue problem
lead to a poor minimum density performance. For example, using $N_{\max}=100$ modes results to a minimum density of 100\%. Therefore, HKS
has to keep 100\% of all correspondence probabilities in the underlying soft correspondence map $S_{[0,1]}$ for enabling to extract still a
perfect matching of 100\%, so that $S_{[0,1]}$ corresponds to a full (dense) matrix.
Therefore, we employ the similarity transform also for the kernel-based methods for a more fair
comparison with the MCR signatures.
\par
\begin{figure}[h]
\centering
\includegraphics[width=0.48\textwidth]{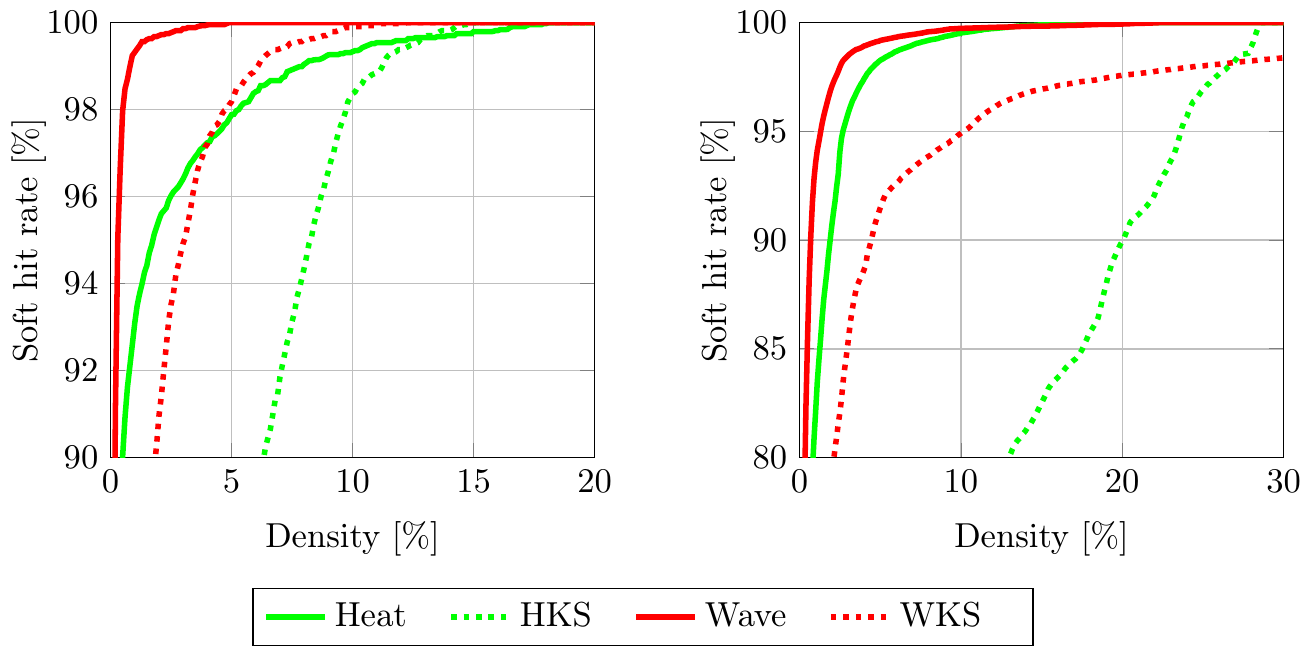}\\
\caption{Results on the dataset wolf (left) and baby (right). 
Comparison of the relation in terms of soft hit rate and the density of $S_{[0,1]}$ between the
  optimised MCR technique and the kernel-based methods for $N_{\max}=100$ modes. The underlying eigenvalues and eigenvectors
  are computed by solving the symmetric eigenvalue problem. We observe that the MCR methods require a much lower number of
  entries in the soft corresponcdence matrix to capture correspondences.}
\label{soft_Wolf_Baby_Comparison}
\end{figure}

We now compare the optimised MCR technique and the kernel-based methods for 100 eigenvalues using the dataset wolf and baby, respectively.
The evaluation of the minimum density highlights the generally better correspondence capability of the MCR method, as shown in
Figure \ref{soft_Wolf_Baby_Comparison} or Table \ref{softmatching}. 
For the wolf and baby dataset the MCR wave and the MCR heat method with 4.86\% and 18.85\% (retained correspondence probabilities for
enabling a perfect matching of 100\%), respectively, outperform all other methods.

\begin{table}[h]
  \caption{Results on the dataset wolf and baby by using the optimised MCR technique and the kernel-based methods for $N_{\max}=100$ modes.
    Comparison of the minimum density 
    of $S_{[0,1]}$ for enabling to extract still a
    perfect matching of 100\%. The corresponding curves are shown in Figure \ref{soft_Wolf_Baby_Comparison}.}
  \label{softmatching}      
    \setlength{\tabcolsep}{10pt}
    \centering\begin{tabular}{|l|ll|}
\hline
        \diagbox{Method}{Dataset}    & Wolf & Baby \\\hline
        MCR heat  & 18.08\% & \textbf{18.85}\% \\[2pt] 
        HKS       & 15.01\% & 28.54\% \\[2pt] 
        MCR wave  & \textbf{4.86\%} & 22.63\% \\[2pt] 
        WKS       & 12.21\% & 62.3\%\\\hline
    \end{tabular}
\end{table}

\paragraph{Evaluation on TOSCA Dataset}

\begin{figure}[t]
\centering
\includegraphics[width=0.38\textwidth]{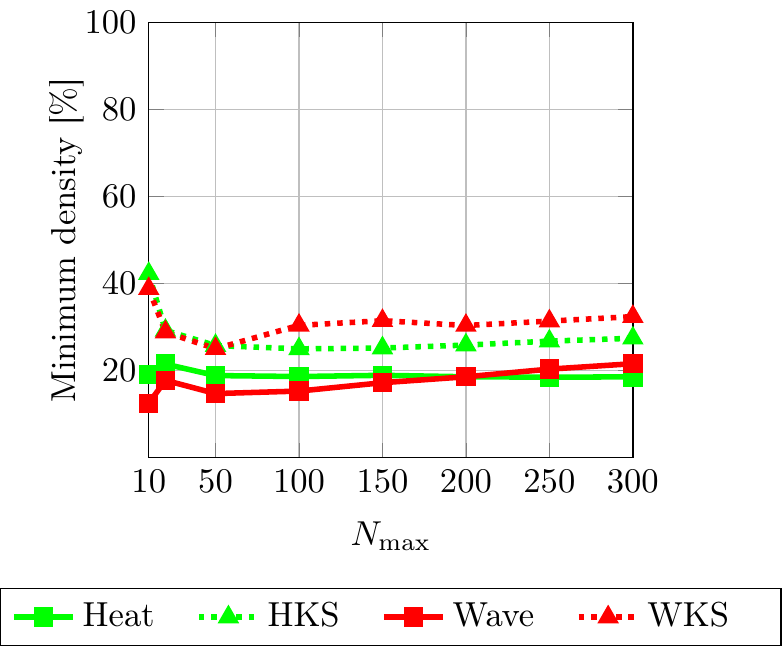}\\
\caption{Results on the entire TOSCA dataset. Comparison of the minimum density between the optimised MCR technique and the
  kernel-based methods for a varying number of modes $N_{\max}$.
  The underlying eigenvalues and eigenvectors are computed by solving the beneficial symmetric eigenvalue problem.
  One clearly observes that the MCR methods rely on a much lower number of entries in the sparsified soft maps.
  Since we average here over the whole TOSCA dataset, this experiment demosntrates the robustness of the MCR approach
  in this respect.
}
\label{Tosca_Comparison_soft}
\end{figure}
\begin{figure}
\centering
\includegraphics[width=0.48\textwidth]{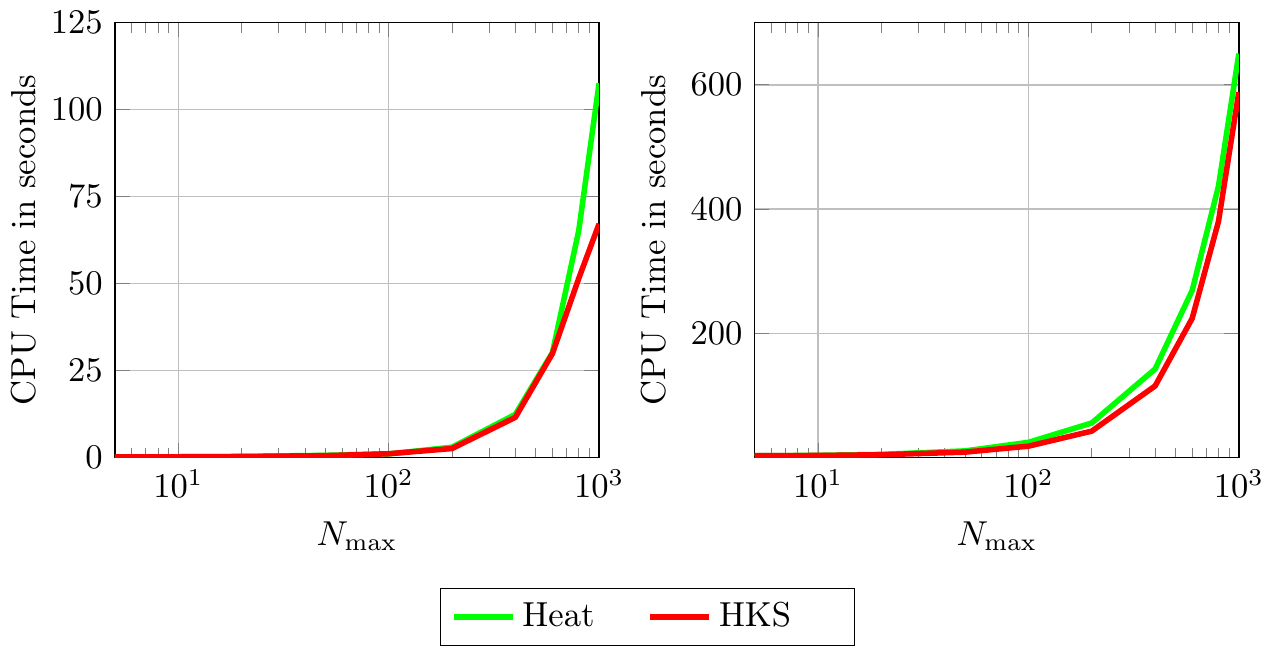}\\
\caption{Results on the datasets wolf (left) and baby (right) by using the geometric heat equation. Comparison of the required computational
  time between the optimised MCR technique and the kernel-based methods for a varying number of modes $N_{\max}\in[5,1000]$. The underlying
  eigenvalues and eigenvectors are computed by solving the symmetric eigenvalue problem. The CPU time of the kernel-based methods is slightly
  faster in particular for $N_{\max}>100$ modes, however the MCR technique is clearly competitive. Using the geometric wave equation generates
  the same costs.}
\label{Comparison_time}
\end{figure}
Also for the whole \\TOSCA dataset, the MCR technique clearly outperforms the kernel-based methods in terms of the minimum density
as shown in Figure \ref{Tosca_Comparison_soft}. In general, the performance of the MCR wave method is the best one among all tested methods,
however from $N_{\max}\approx 200$ modes on the MCR heat method gains in relevance.
The most inferior results are achieved by using WKS, which generates for all used modes $N_{\max}\in[10,300]$ at least 10\% worse density.
\par
As in the previous examination in Section \ref{EvaluationonGeodesicError}, the experiment shows that a small spectrum of
eigenvalues ($N\approx 100$) works well for MCR as well as for the kernel-based methods.
In this regime, the methods MCR wave, MCR heat, HKS and WKS result in required minimum density rates
of 15.29\%, 18.63\%, 24.98\% and 30.39\%, respectively. The latter observation obviously demonstrates
the benefits of using the MCR technique.
\par
In the total, by using a small spectrum the computational effort may be reduced significantly without losing
performance, which is especially important in the case of highly resolved meshes.
\par
In a final experiment, we investigate the required computational time of the applied methods based on the wolf and baby dataset.
As shown in Figure \ref{Comparison_time}, both techniques achieve almost equally fast computational times.
\section{Summary and Conclusion}
We have extended the numerical framework that has been presented in earlier literature,
and in our opinion we have tweaked with this paper the MOR approach very close to its limit
concerning its use in shape correspondence construction. For achieving this we combined
many numerical techniques with benefit that lead to an algorithm that is finally
easy to implement.
\par
Let us stress that our basic approach is nearly free of parameters, one can even conclude
for the remaining few parameters like e.g.\ the number of eigenvalues that we have shown experimentally how
to choose them in applications, so that in practice our approach can be considered as parameter-free, which is
of high practical value.
\par
We think that the use of the MCR technique for particular tasks in shape analysis, like
e.g.\ symmetry or similarity detection, could be a promising subject of future research.
Furthermore, the use of soft correspondences as studied here appears to be promising. As we pointed out and studied here
for the first time in relation to our MOR approach to which it seems to fit quite well, it enables very high matching
results and just needs sparse information for this in the soft correspondence matrix. It is not trivial and
beyond the scope of this work to develop a method that could make full use of the soft correspondence information,
which may result in a quite powerful approach.
\par
Let us also note that as part of the functional maps pipeline mentioned in related work, feature descriptors are used to
compute a coarse correspondence between points. Providing accurate initial correspondences based on feature descriptors
will improve performance and computational time of dense shape correspondence algorithms relying on functional maps, as reported in \cite{MORCC}.
We conjecture that our method may also be useful in this context and we plan to investigate this in future research.


\bibliographystyle{spmpsci}      
\bibliography{references2}   

\end{document}